\documentclass[11pt, reqno]{amsart}
	\usepackage{amssymb, latexsym, amsmath, amsfonts, mathrsfs,bm, amsthm}
	\usepackage{graphicx}
	\usepackage{float}
	\usepackage{tikz}
	\usepackage{tikz-cd}
	\usepackage{enumitem}
	\usepackage[numbers]{natbib}
	\usepackage{hyperref}
	\usepackage{upgreek}
	\makeatletter
\@namedef{subjclassname@2020}{%
  \textup{2020} Mathematics Subject Classification}
\makeatother
	\newtheorem{thm}{Theorem}[section]
	\newtheorem{cor}[thm]{Corollary}
	\newtheorem{lem}[thm]{Lemma}
	\newtheorem{prop}[thm]{Proposition}
 \newtheorem{result}[thm]{Result}
	\theoremstyle{definition}
	\newtheorem{defn}[thm]{Definition}
	\theoremstyle{remark}
	\newtheorem{rem}[thm]{Remark}
 \newtheorem*{rem*}{Remark}
	\numberwithin{equation}{section}
	
	\numberwithin{equation}{section}	
	
\newcommand{\C}{\mathbb{C}}

\newcommand{\Aut}{\mathsf{Aut}(\mathbb{C}^k)}
\newcommand{\seq}[1]{\{#1_n\}}

\newcommand{\h}{\mathsf{H}}
\newcommand{\inv}{\mathsf{i}}
\newcommand{\p}{\mathbf{P}}
\newcommand{\q}{\mathsf{Q}}
\newcommand{\ra}{\rightarrow}
\newcommand{\supnorm}[1]{\|#1\|_{\infty}}
\newcommand{\norm}[1]{\|#1\|}
\newcommand{\abs}[1]{|#1|}
\newcommand{\zball}[1]{B(0;#1)}
\newcommand{\w}[1]{\widetilde{#1}}

\newcommand{\basin}[1]{\Omega_{#1}}
\newcommand{\I}[1]{\textsf{int}(#1)}
\newcommand{\V}[1]{\textsf{vol}(#1)}

\newcommand{\je}{\mathsf{j}}

\newcommand{\li}{\mathsf{L}}
\newcommand{\Rh}{\bm{\rho}}
\newcommand{\Al}{\bm{\alpha}}

\newcommand{\ze}{\mathbf{z}}
\newcommand{\we}{\mathbf{w}}
\newcommand{\ep}{\epsilon}
\newcommand{\co}{(z_1,\hdots,z_k)}
\newcommand{\me}{\mathsf{m}}
\newcommand{\Me}{\mathsf{M}}
\newcommand{\ku}{\mathbf{m}}
\newcommand{\G}{\mathcal{G}}
\newcommand{\A}{\mathsf{A}}
\newcommand{\B}{\mathsf{B}}
\newcommand{\s}{\mathsf{S}}
\newcommand{\un}{\mathbf{u}_n}
\newcommand{\ho}{\mathbf{H}}

\newcommand{\en}{\mathbf{n}}
\begin{document}
\title[Non-autonomous basins]{Uniform Non-autonomous basins of attraction}
	\author{Sayani Bera, Kaushal Verma} 
\address{SB: School of Mathematical and Computational Sciences, Indian Association for the Cultivation of Science, Kolkata-700032, India}
\email{sayanibera2016@gmail.com, mcssb2@iacs.res.in}

\address{KV: Department of Mathematics, Indian Institute of Science, Bengaluru-560012, India}
\email{kverma@iisc.ac.in}

\begin{abstract}
It has been conjectured that every stable manifold arising from a holomorphic automorphism, that acts hyperbolically on a compact invariant set, is biholomorphic to complex Euclidean space. Such stable manifolds are known to be biholomorphic to the basin of a uniformly attracting family of holomorphic maps. It is shown that the basin of a uniformly attracting family of holomorphic maps is biholomorphic to complex Euclidean space and this resolves the conjecture on the biholomorphism type of such stable manifolds affirmatively.
\end{abstract}

\maketitle

\section{Introduction}

\noindent Let $f$ be a holomorphic automorphism of a complex manifold $M$ and suppose that $K \subset M$ is a compact invariant set on which $f$ is uniformly hyperbolic with stable dimension $k$. For $p \in K$, let $\Sigma^p_f$ be the stable manifold of $f$ passing through $p$. Each stable manifold is an immersed complex manifold of dimension $k$ by the stable manifold theorem; in fact, they are diffeomorphic to $\mathbb R^{2k}$. The question of whether $\Sigma^p_f$ is biholomorphic to $\C^k$ for every $p \in K$ was raised by Bedford, and while this is known to be true in several cases (see for instance  \cite{AAM:non-autonomous conjugation}, \cite{BS1}, \cite{JV}, \cite{PS:adaptive}, \cite{RR:holomorphic}), a result of Forn\ae ss--Stens\o nes \cite{FS} shows that the stable manifold $\Sigma^p_f$ is biholomorphic to a domain in $\C^k$ for each $p \in K$. This was done by studying a reformulation of Bedford's question which can be described as follows.

\medskip

Let $\{f_n\}$ be a sequence of holomorphic automorphisms of $\C^k$, $k \ge 2$, such that $f_n(0) = 0$ for all $n \ge 1$. The non-autonomous basin of attraction at the origin is the set defined by
\[
\Omega_{\seq{f}} = \{z \in \C^k: \lim_{n \ra \infty} f_n \circ f_{n-1}\circ\cdots\circ f_1(z) = 0\}.
\]
If there are uniform constants $r > 0$ and $0 < A < B < 1$ such that
\begin{align}\label{e:ub}
A\norm{z} \le \norm{f_n(z)} \le B \norm{z}	
\end{align}
for all $n \ge 1$ and $z \in B(0; r) \subset \C^k$, then $\Omega_{\seq{f}}$ will be said to be a {\it uniform non-autonomous basin of attraction}.

\medskip

\noindent {\it Conjecture}: A uniform non-autonomous basin of attraction is biholomorphic to $\C^k$.

\medskip

A succinct summary of the various partial results obtained thus far that are related to this conjecture can be found in Peters--Smit \cite{PS:adaptive}. In addition, recent work of Forn\ae ss--Wold \cite{FW} shows that $\Omega_{\seq{f}}$ is elliptic and hence an Oka manifold. The following theorem is the main result of the present paper

\begin{thm}\label{t:main result}
A uniform non-autonomous basin of attraction is biholomorphic to $\C^k$.
\end{thm}

Combining \cite[Theorem 3.5]{FS} with the above theorem affirmatively answers Bedford's question, i.e.,

\begin{cor}

Let $f$ be a holomorphic automorphism of a complex manifold $M$ and $K \subset M$ a compact invariant set on which $f$ is uniformly hyperbolic with stable dimension $k \ge 1$. Then, for each $p \in K$, the stable manifold $\Sigma^p_f$ is biholomorphic to $\C^k$.

\end{cor}
Let $k_0\ge 2$ be the smallest positive integer such that $B^{k_0} < A$. This $k_0$ will be fixed once and for all. Also, we will assume that $k_0 \ge 3$, since Theorem \ref{t:main result} is already known to be true for $k_0=2$ by \cite[Theorem 4]{W:FB domains}. Before describing the main steps in the proof, it is useful to recall the work of Rosay--Rudin (\cite{RR:holomorphic}) which shows that the basin of attraction of an attracting fixed point of a holomorphic automorphism of $\C^k$ is biholomorphic to $\C^k$. Lemma $3$ in the Appendix of \cite{RR:holomorphic} shows that for a given integer $m \ge 2$, a holomorphic map of $\C^k$ with an attractive fixed point at the origin can be conjugated to a lower triangular polynomial automorphism of $\C^k$ up to terms of order $m$, by means of polynomial maps of $\C^k$ that are tangent to the identity at the origin. The lower triangular polynomial automorphisms indicated above have the same linear terms at the origin as the given holomorphic automorphism. This idea of creating a conjugacy between the given holomorphic automorphism and a simpler map up to arbitrary orders is effective. It allows for an understanding of the basin of attraction of the holomorphic automorphism in terms of the basin of attraction of a simpler mapping, namely a lower triangular map. This idea has been adapted in the non-autonomous case as well and this forms the basis of several results in this direction.

\medskip

There are two main steps in the proof of the main result. Both build on the framework developed in \cite{AAM:non-autonomous conjugation} and \cite{P:perturbed}. The essential difference in our approach lies not in working with lower triangular maps but compositions of H\'enon maps in $\C^2$ and maps that are close to weak shift like maps in $\C^k$ when $k \ge 3$. 

\medskip

First consider the case $k = 2$. After a suitable normalization of the sequence $\{f_n\}$, the bounds in $(1.1)$ imply the existence of a sequence of uniformly bounded maps $h_n, g_n$ where 
\[
h_n=\textsf{Identity} + \textsf{h.o.t.}
\]
and $g_n$ is the product of two polynomial automorphisms of $\C^2$, each of which is conjugate to a H\'{e}non map, such that 
\[
h_{n+1} = g_n \circ h_n \circ f^{-1}_n
\]
up to order $k_0$ at the origin. In this case, it is known that $\Omega_{f_n} \backsimeq \Omega_{g_n}$. It remains to show that $\Omega_{g_n} \backsimeq \C^2$ and this constitutes the second step. By construction, this family of H\'{e}non maps can be chosen to have uniformly bounded degrees, admit a uniform filtration of $\C^2$ that localizes their dynamics and be uniformly bounded on large enough polydiscs around the origin. The Green's function associated with this family of H\'{e}non maps exists as a non-constant continuous plurisubharmonic function on $\C^2$ and its properties are used to conclude that $\Omega_{g_n} \backsimeq \C^2$. All these points are explained in Section $2$.

\medskip

When $k \ge 3$, while the two steps indicated above remain the same in principle, the details are significantly different. The choice of $g_n$ to which the $f_n$'s must be conjugated to up order $k_0$ is dictated by the following conditions -- first, the equation 
\[
h_{n+1} = g_n \circ h_n \circ f^{-1}_n
\]
which is required to hold up to order $k_0$ at the origin, must admit a solution with 
\[
h_n=\textsf{Identity} + \textsf{h.o.t.}
\]
for all $n$. Second, this choice of $g_n$ must be amenable from the dynamical view point. That is, there must be an analogous Green's function with suitable properties which can be used to show that $\Omega_{g_n}$ is biholomorphic to $\C^k$. Just as H\'{e}non maps play the role of building blocks for the $g_n$'s when $k = 2$, the class of weak shift-like maps are used to identify the $g_n$'s when $k \ge 3$. It turns out that the $g_n$'s can be taken to be the composition of $k$ weak shift-like maps. Perturbing these maps by adding terms of high degree gives rise to a class of maps that are near to them. The basin associated with this sequence of perturbed maps can be shown to be biholomorphic to $\C^k$. Being close to the weak shift-like maps implies that the basins associated with the perturbed maps and the weak shift-like maps are biholomorphic. It follows that $\Omega_{g_n} \backsimeq \C^k$. Sections \ref{s:3} and \ref{s:4} contain the relevant details about this class of maps and this leads to the proof of the main theorem when $k \ge 3$ which is explained in Section \ref{s:5}. Proposition $5.1$ contains the mains inductive step which is required to show the existence of $g_n, h_n$ satisfying 
\[
h_{n+1} = g_n \circ h_n \circ f^{-1}_n
\]
up to order $k_0$ at the origin.

\medskip

The property that a pair of maps have the same finite order jets at a given point is one that will appear frequently, and for this reason, we fix notation as follows: for a $\C^k$-valued holomorphic map $F$ defined near the origin in $\C^k$ and an integer $l \ge 1$, let $[F]_l$ denote the $l$-th order jet of $F$ at the origin. For a pair $F, G: U \ra \C^k$ of holomorphic maps defined on a neighbourhood $U$ of the origin in $\C^k$, the equality of the $l$-th order jets of $F, G$ at the origin will be expressed by writing it as $[F]_l = [G]_l$.

\medskip

To make things clear, Section \ref{s:6} contains a worked-out example that illustrates the main ideas in the inductive step. 

\medskip

\noindent {\it Acknowledgements}: We wish to thank the referees for providing several excellent suggestions and comments. They have resulted not only in making the exposition clearer but have also contributed in simplifying the proofs of some intermediate steps leading to the proof of the main theorem. This is specially true of Section \ref{s:5}, much of which stands revised and simplified thanks to a rather pertinent suggestion. We would also like to thank Gautam Bharali for several useful discussions.

\section{A Computation in \texorpdfstring{$\mathbb C^2$}{C2}}\label{s:2}

As indicated in \cite[Remark 7]{P:perturbed}, by working with the $QR$-factorization of $Df_n(0)$, $n \ge 1$, the linear part of every $f_n$ can be assumed to be a lower triangular matrix. In particular,
\[Df_n(0)=\begin{pmatrix} a_n &0\\b_n &c_n \end{pmatrix} \text{ and }Df_n^{-1}(0)=\begin{pmatrix} a_n^{-1} &0\\b'_n &c_n^{-1} \end{pmatrix} \]
for every $n \ge 1$. By the bounds assumed in (\ref{e:ub}), it follows that
\begin{align}\label{e:inverse ub}
B^{-1}\norm{z} \le \norm{f_n^{-1}(z)} \le A^{-1}\norm{z}.
\end{align}
on $B(0;Ar)$. The Cauchy estimates applied to the components of $f_n$ along coordinate discs of radius $r$ in $B(0;r)$ give
\[
0\le |a_n|\le B, \; 0 \le |b_n|\le B \text{ and } 0 \le |c_n|\le B.
\]
Similarly, the above bounds on the growth of $f_n^{-1}$ on the ball $B(0;Ar)$ give 
\[
0\le |a_n^{-1}|\le A^{-1}, \; 0 \le |b'_n|\le A^{-1} \text{ and } 0 \le |c^{-1}_n|\le A^{-1}.
\]
Hence
\begin{itemize}
	\item $A\le |a_n|\le B$, $A\le |c_n|\le B$ and $0 \le \vert b_n \vert \le B$,
	\item $|a_n||c_n|^{-k}\ge c_0>1$ and  $|c_n||a_n|^{-k}\ge c_0>1$ for $k \ge k_0$.
\end{itemize}
for every $n \ge 1$. 

\medskip Our main goal here is to identify a sequence of maps $\{g_n\}$ such that the sequences $\{f_n\}$ and $\{g_n\}$ are {\it non-autonomously conjugated up to order $k_0$}.  It will turn out that the $g_n$'s can be chosen to be H\'{e}non maps. To recall the definition of {\it non-autonomous conjugation up to a given order}, we first note the following: a sequence $\{\varrho_n\}$ of holomorphic self-maps of $\C^k$ is said to be {\it uniformly bounded at the origin} if the restrictions of $\varrho_n$ to some uniform ball $B(0;\delta)$, $\delta > 0$, are uniformly bounded. 

\begin{defn}
 A pair of sequences $\{f_n\}, \{g_n\}$, both of which are uniformly bounded automorphisms of $\C^k$, $k \ge 2$, are said to be 
 {\it non-autonomously conjugated up to order $m$}, $m \ge 1$,  if there exists a sequence $\seq{h}$ of endomorphisms of $\C^k$, $k \ge 2$, of the form 
\[h_n=\textsf{Identity} + \textsf{h.o.t.}\] 
near the origin such that $\{h_n\}$ is uniformly bounded at the origin and
each box in the diagram
\[
\begin{tikzcd}
 \C^2 \arrow [r,"f_1"]\arrow[d, "h_1"] & \C^2  \arrow[r, "f_2"] \arrow[d,"h_2"] & \C^2 \arrow[r, "f_3"] \arrow[d,"h_3"]& \cdots \\
 \C^2 \arrow [r,"g_1"] & \C^2  \arrow[r, "g_2"] &\C^2  \arrow[r, "g_3"] & \cdots    
 \end{tikzcd}
\] 
commutes up to order $m$  at the origin; that is, $[h_{n+1}]_m = [g_n \circ h_n \circ f^{-1}_n]_m$.    
\end{defn}

\begin{prop}\label{p:step 1}
There exist sequences of polynomials $\seq{p}$ and $\seq{q}$ of degree $k_0$ in one variable with 
\[p_n(0)=q_n(0)=p_n'(0)=0 \text{ and }q_n'(0)=b_n \] 
 for every $n \ge 1$ such that the sequence $\seq{g}$ defined by 
\begin{align}\label{e:rseq}
g_n(x,y)= \Big(a_n x+p_n(y),c_n y+q_n\big(x+a_n^{-1}p_n(y)\big) \Big) 
\end{align}
is non-autonomously conjugated to $\seq{f}$ up to order $k_0$.
\end{prop}

\begin{proof}
To begin with, for every $2 \le m \le k_0$, let 
\[
\mathrm{I}_{m}=\{(i,j):i,j \in \mathbb{N} \cup\{0\} \text{ and }2 \le i+j\le m\}
\]
be the set of all multi-indices up to degree $m$ in two variables and 
$N_m=\sharp(\mathrm{I}_m)$ the cardinality of $\mathrm{I}_m$. Let
 \begin{align*}
	 f_n(x,y)&=\Big(a_nx+\sum_{(i,j)\in \mathrm{I}_{k_0}}\mu_{n,i,j}^1x^iy^j,b_nx+c_ny+\sum_{(i,j)\in \mathrm{I}_{k_0}}\mu^2_{n,i,j}x^iy^j\Big)+\textsf{h.o.t.}
  \end{align*}
and 
\begin{align*}
	 f_n^{-1}(x,y)&=\Big(a_n^{-1}x+\sum_{(i,j)\in \mathrm{I}_{k_0}}\theta_{n,i,j}^1x^iy^j,b'_nx+c_n^{-1}y+\sum_{(i,j)\in \mathrm{I}_{k_0}}\theta^2_{n,i,j}x^iy^j\Big)+\textsf{h.o.t.}
  \end{align*}
  From (\ref{e:ub}) and (\ref{e:inverse ub}), there exists $M>0$ large enough such that for every $l=1,2$, $n \ge 1$ and $(i,j) \in \mathrm{I}_{k_0}$
  \begin{align}\label{e:coefficients} 
  0\le |\mu_{n,i,j}^l|<M \text{ and } 0 \le |\theta_{n,i,j}^l|<M
  \end{align}
  Further, for $z \in \C$ let
  \begin{align}\label{e:p and q} 
  p_n(z)=\sum_{i=2}^{k_0}\alpha^1_{n,0,i}z^i \text{ and } q_n(z)=b_n z+\sum_{i=2}^{k_0}\alpha^2_{n,i,0}z^i.
  \end{align}Then
  \begin{align*}
	{g}_n(x,y)&=\Big(a_nx+\sum_{j=2}^{k_0} \alpha^1_{n,0,j}y^j,c_ny+b_n\big(x+a_n^{-1}p_n(y)\big)+\sum_{i=2}^{k_0}\alpha^2_{n,i,0}\big(x+a_n^{-1}p_n(y)\big)^i\Big).
\end{align*}
Let $\alpha^1_{n,i,j}=0$ for $i \ge 1$ and $\alpha^2_{n,i,j}=0$ for $j \ge 1$. Also, let $\seq{{h}}$ be maps of the form
\begin{align*}
	\pi_1\circ {h}_n(x,y)=x+\sum_{(i,j)\in \mathrm{I}_{k_0}}\rho^1_{n,i,j}x^iy^j\text{ and }
	\pi_2\circ {h}_n(x,y)=y+\sum_{(i,j)\in \mathrm{I}_{k_0}}\rho^2_{n,i,j}x^iy^j.
\end{align*}
where the projections $\pi_i$ capture the first and second components of $h_n$.

\medskip

\noindent For every $n \ge 1$, $2 \le m \le k_0$ and $l=1,2$, consider the following tuples 
\[ \Al_{n,m}^l=(\alpha^l_{n,i,j}),\Rh^l_{n,m}=(\rho^l_{n,i,j})\in \C^{m+1} \text{ where } m=i+j.\] 
Note that $\Al^1_{n,m}=(0,\hdots,0,\alpha^1_{n,0,m})$ and $\Al^2_{n,m}=(\alpha^2_{n,m,0},0,\hdots,0).$ Also, let
\[ \Rh_{n,m}=(\Rh^1_{n,2},\hdots,\Rh^1_{n,m}, \Rh_{n,2}^2,\hdots, \Rh^2_{n,m}),\;  \Al_{n,m}=(\Al^1_{n,2},\hdots,\Al^1_{n,k}, \Al_{n,2}^2,\hdots, \Al^2_{n,m})\in \C^{N_m}.\] 
Hence, the maps $\seq{{h}}$ and $\seq{{g}}$ can be uniquely identified with the sequences $\{\Rh_{n,k_0}\}$ and $\{\Al_{n,k_0}\}$ respectively. We will construct these sequences inductively by using the relation
\begin{align}\label{e:inductive}
[{h}_{n+1}]_{k_0}=[{g}_n \circ {h}_n \circ f_n^{-1}]_{k_0}.
\end{align} 

 To complete the proof, we need to ensure that there exist bounded sequences $\{\Rh_{n,k_0}\}$ and $\{\Al_{n,k_0}\}$ such that (\ref{e:inductive}) holds.   

\medskip

For $i+j=2$, we need to compare the degree two terms in (\ref{e:inductive}). To do this, first note that the degree two terms of $h_n \circ f_n^{-1}$ contain contributions by both the degree one and degree two terms of $f_n^{-1}$ and $h_n$. Similarly, the degree two terms of $g_n \circ h_n \circ f_n^{-1}$ contain contributions by both the degree one and degree two terms of $h_n \circ f_n^{-1}$ and $g_n$. Now observe that 

\[[\pi_1 \circ h_n \circ f_n^{-1}(x,y)]_2=a_n^{-1}x+\sum_{i+j=2} \theta_{n,i,j}^1 x^i y^j+\sum_{i+j=2} \rho_{n,i,j}^1 (a_n^{-1} x)^i (c_n^{-1}y+b_n' x)^j \] and
\[[\pi_2 \circ h_n \circ f_n^{-1}(x,y)]_2=c_n^{-1}y+b_n' x+\sum_{i+j=2} \theta_{n,i,j}^2 x^i y^j+\sum_{i+j=2} \rho_{n,i,j}^2 (a_n^{-1} x)^i (c_n^{-1}y+b_n' x)^j. \]
Further 
\begin{align} \label{e:C2 degree 2_1} 
\nonumber [\pi_1 \circ  g_n\circ h_n \circ f_n^{-1}(x,y)]_2= x&+a_n\sum_{i+j=2} \theta_{n,i,j}^1 x^i y^j+a_n\sum_{i+j=2} \rho_{n,i,j}^1 (a_n^{-1} x)^i (c_n^{-1}y+b_n' x)^j\\  &+\alpha_{n,0,2}^1(c_n^{-1}y+b_n'x)^2
\end{align}
and 
\begin{align} \label{e:C2 degree 2_2}
\nonumber [\pi_2 \circ  g_n\circ h_n \circ f_n^{-1}(x,y)]_2= y&+c_n\sum_{i+j=2} \theta_{n,i,j}^2 x^i y^j+c_n\sum_{i+j=2} \rho_{n,i,j}^2 (a_n^{-1} x)^i (c_n^{-1}y+b_n' x)^j\\ 
\nonumber &+b_n\sum_{i+j=2} \theta_{n,i,j}^1 x^i y^j+b_n\sum_{i+j=2} \rho_{n,i,j}^1 (a_n^{-1} x)^i (c_n^{-1}y+b_n' x)^j\\ &+b_na_n^{-1} \alpha_{n,0,2}^1 (c_n^{-1}y+b_n'x)^2+\alpha_{n,2,0}^2 (a_n^{-1}x)^2
\end{align}
Thus, comparing the terms of (\ref{e:C2 degree 2_1}) and (\ref{e:C2 degree 2_2}) with (\ref{e:inductive}), we get the relations:
\begin{align*}
\rho_{n+1,0,2}^1&=a_nc_n^{-2}\rho_{n,0,2}^1+c_n^{-2}\alpha_{n,0,2}^1+\li_{n,0,2}^1,\\
\rho_{n+1,1,1}^1&=c_n^{-1}\rho_{n,1,1}^1+\li_{n,1,1}^1(\rho_{n,0,2}^1,\alpha_{n,0,2}^1),\\
\rho_{n+1,2,0}^1&=a_n^{-1}\rho_{n,0,2}^1+\li_{n,0,2}^1(\rho_{n,0,2}^1,\rho_{n,1,1}^1,\alpha_{n,0,2}^1),\\
\rho_{n+1,0,2}^2&=c_n^{-1}\rho_{n,0,2}^2+\li^2_{n,0,2} (\rho_{n,0,2}^1,\alpha_{n,0,2}^1),\\
\rho_{n+1,1,1}^2&=a_n^{-1}\rho_{n,1,1}^2+\li^2_{n,1,1}(\rho_{n,0,2}^1,\rho_{n,1,1}^1,\rho_{n,0,2}^2,\alpha_{n,0,2}^1),\\
\rho_{n+1,2,0}^2&=a_n^{-2}c_n\rho_{n,2,0}^2+\li^2_{n,2,0} (\rho_{n,0,2}^1,\rho_{n,1,1}^1 , \rho_{n,2,0}^1,\rho_{n,0,2}^2,\rho^2_{n,1,1},\alpha_{n,0,2}^1)+a_n^{-2}\alpha_{n,2,0}^2, 
\end{align*}
where $\li^l_{n,i,j}$ is a linear combination of the $\rho^l_{n,i,j}$'s with coefficients depending on those of $f_n^{-1}$ up to degree 2. By (\ref{e:coefficients}), these degree $2$ coefficients of $f_n^{-1}$ are uniformly bounded.  Let $\{\alpha^1_{n,0,2}\}$ be the bounded sequence such that $\rho_{n,0,2}^1=0$ for every $n \ge 1.$ To obtain a bounded solution for the rest of the terms $\{\Rh_{n,2}\}$, except $\rho_{n,2,0}^2$, we will need the following lemma from \cite[Lemma 12]{P:perturbed}, and it will be repeatedly used to move from one equation to the next in the list above.

\begin{lem}\label{l:linear bound}
Let $\seq{A}$ be a sequence of bounded strictly expanding affine maps of $\C$ of the form $A_n(z)=\beta_n z+\gamma_n$ where $1<c<|\beta_n|<C,\; 0\le  |\gamma_n|<C$ for uniform constants $1 < c < C$.
Then there exists a bounded orbit $\{z_n\}$ of $\seq{A}$.
\end{lem}
 
 As $|c_n^{-1}|>B^{-1}$, there exists an appropriate $\rho_{0,1,1}^1$ such that $\{\rho^1_{n,1,1}\}$ is a bounded orbit, in particular, it is a bounded sequence. Once the sequence $\{\rho^1_{n,1,1}\}$ is fixed, by appealing to the above lemma again, it is possible to choose $\rho_{0,2,0}^1$ such that $\{\rho^1_{n,2,0}\}$ is bounded. Similarly, there exist bounded sequences $\{\rho_{n,0,2}^2\}$, $\{\rho_{n,1,1}^2\}$. Now for the last equation in the list, we apply the same idea that was used for $\rho_{n,0,2}^1$, i.e., choose the bounded sequence $\{\alpha_{n,2,0}^2\}$ appropriately such that $\rho_{n,2,0}^2=0$ for every $n \ge 1$.

\medskip 

Note that for a fixed $m$, $2<m=i+j\le k_0$, the degree $m$ terms in the expansion of $h_{n+1}$ are made up of all the terms up to degree $m-1$ in the expansions of $f_n^{-1}$, $h_n$ and $g_n$, and a linear combination of degree $m$ terms in the expansions of $f_n^{-1}$, $h_n$ and $g_n$. Let $\mathcal{A}_{n,m-1}^l(x,y)$, $\mathcal{R}_{n,m-1}^l(x,y)$ and $\Theta_{n,m-1}^l(x,y)$ denote the sum of degree two to degree $m-1$ terms in the expansions of  $\pi_l \circ g_n$, $\pi_l \circ h_n$ and $\pi_l \circ f_n^{-1}$, for $l=1,2$, respectively. Then 
{\small \begin{align*}
[g_n(x,y)]_m=&\Big(a_n x+\mathcal{A}_{n,m-1}^1(x,y)+ \alpha_{n,0,m}y^m,\; b_nx+c_ny+\mathcal{A}_{n,m-1}^2(x,y)+\alpha_{n,m,0}^2 x^m\Big)\\
[h_n(x,y)]_m=&\Big(x+\mathcal{R}_{n,m-1}^1(x,y)+\sum_{i+j=m} \rho_{n,i,j}^1x^iy^j,\; y+\mathcal{R}_{n,m-1}^2(x,y)+\sum_{i+j=m} \rho_{n,i,j}^2 x^iy^j\Big)\\
[f_n^{-1}(x,y)]_m=&\Big(\Theta_{n,m-1}^1(x,y)+\sum_{i+j=m} \theta_{n,i,j}^1x^iy^j\; ,+\Theta_{n,m-1}^2(x,y)+\sum_{i+j=m} \theta_{n,i,j}^2x^iy^j\Big)+\\&(a_n^{-1}x,\; b'_nx+c_n^{-1}y)
\end{align*}
\noindent So if we assume that there exist bounded sequences $\{\Rh_{n,m-1}\}, \{\Al_{n,m-1}\}$, then $\mathcal{A}_{n,m-1}^l(x,y)$ and $\mathcal{R}_{n,m-1}^l(x,y)$ are well defined. Also $\Theta_{n,m-1}^l(x,y)$ is always defined and its coefficients are bounded by (\ref{e:coefficients}). Hence, from (\ref{e:inductive}) the sequences $\{\Rh_{n,m}^1\}$, $\{\Rh_{n,m}^2\}$, $\{\Al_{n,m}^1\}$, $\{\Al_{n,m}^2\}$ satisfy the following relations: }
{\small \begin{align*}
\rho_{n+1,0,m}^1&=a_nc_n^{-m}\rho_{n,0,m}^1+c_n^{-m}\alpha_{n,0,m}^1+\li_{n,0,m}^1(\Rh_{n,m-1},\Al_{n,m-1}),\\
\rho_{n+1,i,j}^1&=a_n^{-i+1}c_n^{-j}\rho_{n,i,j}^1+\li_{n,i,j}^1(\Rh_{n,m-1},\Al_{n,m-1},\Al^1_{n,m}, \rho_{n,i-1,j+1}^1,\hdots,\rho_{n,0,m}^1), 1\le i\le m\\
\rho_{n+1,0,m}^2&=c_n^{-m+1}\rho_{n,m,0}^2+\li^2_{n,0,m} (\Rh_{n,m-1},\Al_{n,m-1}, \Rh_{n,m}^1,\Al^1_{n,m})\\
\rho_{n+1,i,j}^2&=c_n^{-j+1}a_n^{-i}\rho_{n,i,j}^2+\li^2_{n,i,j}(\Rh_{n,m-1},\Al_{n,m-1}, \Rh_{n,m}^1,\Al^1_{n,m},\rho_{n,i-1,j+1}^2,\hdots,\rho_{n,0,m}^2), 1\le i\le m-1\\
\rho_{n+1,m,0}^2&=a_n^{-m}c_n\rho_{n,m,0}^2+\li^2_{n,m,0} (\Rh_{n,m-1},\Al_{n,m-1}, \Rh_{n,m}^1,\Al^1_{n,m},\rho_{n,m-1,1}^2,\hdots,\rho_{n,0,m}^2)+a_n^{-m}\alpha_{n,m,0}^2, 
\end{align*}}
where $\li^l_{n,i,j}$ is a polynomial of bounded degree in $\rho^l_{n,i,j},\Rh_{n,i}, \Al_{n,i}$ depending on the coefficients of $f_n^{-1}$, up to degree $m$ (which are uniformly bounded). Thus, it is possible to inductively define bounded sequences $\{(\Rh^1_{n,m},\Rh^2_{n,m}, \Al^1_{n,m},\Al^2_{n,m})\}$ for every $1 \le m \le (k_0-1)$ such that $\rho_{n,0,m}^1=\rho_{n,m,0}^2=0$ for every $n \ge 1$ and $2 \le m \le (k_0-1)$. 

\medskip

Now as $\vert a_nc_n^{-k_0} \vert, \vert c_na_n^{-k_0} \vert > 1$ and
$\{(\Rh_{n,k_0-1}, \Al_{n,k_0-1})\}$ are bounded sequences, there exists a bounded solution $\{\Rh_{n,k_0}\}$ such that $\alpha_{n,0,k_0}^1=\alpha_{n,k_0,0}^2=1$ and $\alpha_{n,i,j}^1=\alpha_{n,i,j}^2=0$ otherwise with $i+j=k_0$.

\medskip\noindent
 Hence, the sequences $p_n$ and $q_n$, which are monic polynomials of degree $k_0$ by (\ref{e:p and q}), are uniformly bounded at the origin with $p_n'(0)=0$, $q_n'(0)=b_n$.  By construction, the sequence $g_n(x,y)= (a_n x+p_n(y),c_n y+q_n(x+a_n^{-1}p_n(y))) $ satisfies (\ref{e:inductive}) and this completes the proof.
 \end{proof}

\begin{rem*} Finally, a word of caution. Proposition \ref{p:step 1} does not show that any automorphism of $\C^2$ can be expressed as the composition of two H\'{e}non maps of the form (\ref{e:hm}). In fact, it only says that a finite order jet at the origin of a sequence of automorphisms $\{f_n\}$ satisfying (\ref{e:ub}) can be realised with the help of an appropriate sequence of H\'{e}non maps $\{g_n\}$ and a sequence of bounded polynomial maps $\{h_n\}$, with $Dh_n(0)=\textsf{Id}$, both of degrees at most $k_0$.
\end{rem*}

 \noindent To move ahead, a paraphrasing of Theorem A.1 and Lemma 5.2 from \cite{AAM:non-autonomous conjugation} gives 
 
 \begin{result}
Let $\{f_n\}$, $\{g_n\}$ be sequences of bounded holomorphic functions from a fixed neighbourhood of the origin in $\C^k$ to $\C^k$, $k \ge 2$, that are uniformly bounded at the origin. Assume that 
on the ball $B(0;r)$ both are uniformly attracting at the origin with coefficients $A$ and $B$ (as in \ref{e:ub}), i.e.,
   \begin{align*}\label{e:ub}
A\norm{z} \le \norm{f_n(z)} \le B \norm{z}	\text{ and } A\norm{z} \le \norm{g_n(z)} \le B \norm{z}
\end{align*}
for $z \in B(0;r)$. If $\seq{f}$ and $\seq{g}$ are non-autonomously conjugate up to order $k_0$, where $B^{k_0}<A<B<1$, then the  abstract basins of attraction of the sequences $\seq{f}$ and $\seq{g}$ are biholomorphic.
 \end{result}
 
 Note that if the sequences $\seq{f}$ and $\seq{g}$ in the above result are germs of automorphisms of $\C^k, k \ge 2$ then the abstract basin of $\seq{f}$ and $\seq{g}$ are known to be naturally biholomorphic to their respective basins of attraction at the origin (See \cite[Remark A.4]{AAM:non-autonomous conjugation}). Hence we have the following
 
 \begin{cor}
 Let $\seq{f}$ and $\seq{g}$ be as in Proposition \ref{p:step 1}. Then $\Omega_{\seq{f}} \backsimeq \Omega_{\seq{g}}$.
 \end{cor}

We will complete the proof of Theorem \ref{t:main result} when $k = 2$ by an application of Corollary 7.4 from \cite{B:semigroup} for which we first recall the following definition:

\medskip 

For $R > 0$, consider the filtration of $\C^2$ given by
\[
V_R = \{\vert x \vert, \vert y \vert < R\}, V^+_R = \{ \vert y \vert \ge \max\{\vert x \vert, R\}\}, V^-_R = \{\vert x \vert \ge \max\{\vert y \vert, R \}\}.
\]
A sequence of generalized H\'{e}non maps $\{\h_n\}$, (that is, $\h_n$ is a finite composition of maps of form 
\begin{align}\label{e:hm}
\h(x,y)=(y,\delta x+P(y))	
\end{align}
where $\delta\neq 0$ and $P$ is a polynomial in one variable of degree at least $2$), is said to satisfy the \emph{uniform filtration and bound condition} if 
\begin{enumerate}[leftmargin=14pt]
\item[(i)] $\seq{\h}$ admits a uniform filtration radius $R_{\seq{\h}}>1$ (sufficiently large) such that for every $R>R_{\seq{\h}}$
\begin{itemize}
	\item[(a)] $\displaystyle \overline{\h_n(V_R^+)} \subset V_R^+ \text{ and } \overline{\h_n^{-1}(V_R^-}) \subset V_R^-.$
	\item[(b)]there exists a sequence positive real numbers $\seq{R}$ diverging to infinity, with $R_0=R$, satisfying
	$\displaystyle V_{R_n} \cap \h(n)(V_R^+)=\emptyset \text{ and }V_{R_n} \cap \h(n)^{-1}(V_R^-)=\emptyset.$
	\item[(c)] there exist uniform positive constants $0<\me<1<\Me$ such that 
	\[\me \abs{y}^{d_n}<\norm{\h_n(x,y)}=\abs{\pi_2 \circ \h_n(x,y)}<\Me \abs{y}^{d_n} \text{ on }V_R^+\]
and
\[\me \abs{x}^{d_n}<\norm{\h_n^{-1}(x,y)}=\abs{\pi_1 \circ \h_n^{-1}(x,y)}<\Me \abs{x}^{d_n} \text{ on }V_R^-.\]
 where $d_n$ is the degree of $\h_n.$ 
 \end{itemize}
 \item[(ii)] For every $R\ge R_{\seq{\h}}$, there exists a uniform constant $B_R=\max\{\norm{\h_n(z)}: z\in V_R\} < \infty.$
 \end{enumerate}

\smallskip

Now Corollary 7.4 in \cite{B:semigroup} is the following 
\begin{result}\label{c:BC}
Let $\{\h_n\}$ be a sequence of H\'{e}non maps of form (\ref{e:hm})
satisfying the uniform filtration and bound condition and is uniformly attracting on a neighbourhood of origin, i.e., satisfying (\ref{e:ub}). Then the basin of attraction of the sequence $\{\h_n\}$ at the origin is biholomorphic to $\C^2$.
\end{result}
\begin{rem}
The above is stated as a corollary to Theorem 1.1 in \cite{B:semigroup} in light of Remarks 6.8 and 6.9 mentioned therein and hence we will not prove it here. However the class of weak shift-like maps and perturbed weak shift-like maps in $\C^k, k \ge 3$, that are introduced in Section \ref{s:3} below, coincides with the class of H\'{e}non maps in $\C^2$ when the number of variables is restricted to two. With this observation, that will become evident in Section \ref{s:3}, Proposition \ref{p:na p-weak-shift basin} in Section \ref{s:4} is a generalization of Result \ref{c:BC} to higher dimensions and the proof works for $k=2$ as well. Also note Remark \ref{r:Henon} in this context.
\end{rem}
\begin{proof}[Proof of Theorem \ref{t:main result} in $\C^2$]Note that
\begin{equation*}
g_n(x,y) =\Big(a_n x+p_n(y),c_n y+q_n\big(x+a_n^{-1}p_n(y)\big)\Big)
\end{equation*}
can be factored as $g_n(x, y)={\h}_{2n} \circ {\h}_{2n-1}(x,y)$ where
\[{\h}_{2n}=\big(y, c_n x+q_n(a_n^{-1}y)\big)\text{ and }{\h}_{2n-1}(x,y)=\big(y,a_n x+p_n(y)\big)\]
for every $n \ge 1$. Thus $g_n(x,y)$ is a H\'{e}non map of the above form (\ref{e:hm}) and the coefficients of the monic polynomials $p_n(x)$ and $q_n(x)$ are bounded by construction. Also $A<|a_n|,|c_n|<B$. Hence, the sequence $\{\h_n\}$ satisfies the uniform filtration and bound condition stated above. Also, the origin is a uniformly attracting fixed point of the sequence $\{g_n\}$. Thus, the proof follows by Result \ref{c:BC}.
\end{proof}

\section{Dynamics of perturbed weak shift-like maps}\label{s:3}

 Our goal here and in the consecutive sections will be to identify non-autonomous sequences of a certain class of polynomial automorphisms of $\C^k, k\ge 3$, for which the above Result \ref{c:BC} and the idea of non-autonomous conjugacy can be achieved. A natural generalizaton of H\'{e}non maps to higher dimensions, for which the analogous potential theoretic properties (based on the existence of a filtration of the ambient space) of iterated dynamical systems are known to extend to, are the {\it shift-like maps}. They were introduced and first studied by Bedford--Pambuccian in \cite{BP:shift maps}. However, the form of shift-like maps does not allow the full range of desired properties to be achieved. A suitable perturbation is necessary and this is described in what follows.

\medskip 

Recall that H\'{e}non maps extend as birational maps to $\mathbb{P}^2$ with non-intersecting indeterminacies. The analogue class of maps in $\C^k, k \ge 3$, are the {\it regular and weakly regular maps} that were introduced and studied in \cite{S:regular} and  \cite{SG:regular}. We will first study some pluripotential theoretic properties of the following class of weakly regular polynomial automorphisms in $\C^k, k \ge 3$, that consists of maps of the form
\[S\co=\big(z_2,\hdots,z_{k-1},z_k+Q_{d-1}(z_2),az_1+p(z_2,\hdots,z_k)+\ho_d(z_2,\hdots, z_k)\big)\]
where $Q_{d}(z)=z^d$, $a \not= 0$,
\[
\ho_d(z_1,\hdots,z_{k-1})=\sum_{i=1}^{k-1}z_i^d
\]
and $d= \max\{\text{deg}(p)+2,3\}.$ Note that 
\[
S^{-1}\co=\Big(a^{-1}\big(z_k+(p+\ho_d)(z_1,\hdots, z_{k-1}-Q_{d-1}(z_1)\big),z_1,\hdots,z_{k},z_{k-1}-Q_{d-1}(z_1)\Big).\]
Note that the map $S$ above when restricted $\C^2$, gives a simple H\'{e}non map of the form 
\[\h(x,y)=(y, ax+p(y)+y^d), \; a \neq 0.\]
\begin{lem}\label{l:wreg}
$S$ and $S^{-1}$ are weakly regular and algebraically stable.	
\end{lem}

\begin{proof}
Denote by $\overline{S}$ and $\overline{S^{-1}}$ the birational extensions of $S$ and $S^{-1}$ respectively to $\mathbb P^k$. Let $\mathcal{H} \subset \C^{k-1}$ be the affine variety defined by the vanishing of the homogeneous polynomial $\ho_d$. The indeterminacy loci of $S$ and $S^{-1}$ are
\begin{align*}
I^+&=\{[0:\ze:0]: \ze \in \mathcal{H}\setminus \{0\} \} \cup \{[1:\ze:0]: \ze \in \mathcal{H}\},\\
I^-&=\{[0:\ze:0]: \ze \in \C^{k-1}\setminus \{0\} \} 
\end{align*}
Note that the highest degree term in the $n$-fold composition $S^n$ is given by the homogeneous polynomial $\left( \ho_d(\ze) \right)^{d^{n-1}}$ where $\ze \in \C^{k-1}$. Similarly, the highest degree term for $S^{-n}$ is given by large enough powers of $z_1$ in the first coordinate of $S^{-n}$. Hence, if $I_n^{\pm}$ denote the indeterminacy loci of $S^{\pm n}$ respectively, it follows that $I_n^+=I^+$ and $I_n^-=I^-$ for every $n \ge 1$. Note that 
\[X_1^+=\overline{S} \left([z_1:\cdots:z_{k-1}:z_k:0]\setminus I^+ \right)=[0:\cdots:0:1:0] \in I^- \setminus I^+\]
Thus, inductively for $n>1$,
\[X_n^+=\overline{S} \left(X_{n-1}^+\setminus I^+ \right)=[0:\cdots:0:1:0] \in I^- \setminus I^+\]
and 
\[X_n^-=\overline{S^{-1}} \left(X_{n-1}^-\setminus I^- \right)=[1:0\cdots:0:0] \in I^+ \setminus I^-.\]
This completes the proof.
\end{proof}

\begin{rem}\label{r:na_wreg}
Consider a family of such maps: 
\begin{align}\label{e:tampered weak-shifts}
S_n\co=\big(z_2,\hdots,z_k+Q_{d_n-1}(z_2),a_nz_1+p_n(z_2,\hdots,z_k)+\ho_{d_n}(z_2,\hdots, z_k)\big)	
\end{align}
where $a_n \neq 0$ and $d_n\ge \text{degree of }(a_n z_1+ p_n(z_2,\hdots,z_k))+2.$ Let $S(n) = S_n \circ \cdots \circ S_1$ and $S^{-1}(n) = S^{-1}_n \circ \cdots \circ S^{-1}_1$. By a similar argument as in the proof of Lemma \ref{l:wreg}, it follows that $S(n), S^{-1}(n)$ are both weakly regular and satisfy
\[X_n^+=\overline{S_n} \left(X_{n-1}^+\setminus I^+\right )=[0:\cdots:0:1:0] \in I^- \setminus I^+\]
and 
\[X_n^-=\overline{S^{-1}_n}\left (X_{n-1}^-\setminus I^-\right )=[1:0\cdots:0:0] \in I^+ \setminus I^-.\]
\end{rem}

\medskip

Since $X_n^+=X^+=[0:\cdots:0:1:0]$ and $X_n^-=X^-=[1:0:\cdots:0:0]$ in the above cases (i.e., both iterative and non-autonomous), consider the following subsets of $\C^k$:
\[U^+=\{ z\in \C^k: \overline{S(n)}([z:1])\to X^+\} \text{ and }U^-=\{ z\in \C^k: \overline{S^{-1}(n)}([z:1])\to X^-\}.\]
Also, let $\mathcal{K}^\pm=\C^k \setminus U^\pm$ and $K^\pm=\{z \in \C^k: S^{\pm}(n)(z) \text{is bounded}\}$. In the iterative case (i.e., $S_n=S$), Theorem 2.2 in \cite{SG:regular} shows that there exist continuous plurisubharmonic Green's functions $G^{\pm}$ on $\C^k$ such that $\mathcal{K}^\pm = \{G^\pm(z)=0\}$. We will first generalize this fact in the setting of appropriate non-autonomous sequences $\seq{S}$ in $\C^k$. Consider the neighbourhoods of $X^\pm$ defined by
\begin{align*}
V_{R}^k&=\{z \in \C^k:|z_k| \ge  \max\{(k-1)|z_1|,\hdots,(k-1)|z_{k-1}|,R\}\}\\ 
V_{R}^-&=\{z \in \C^k:|z_1| \ge  \max\{(k-1)|z_2|,\hdots, (k-1)|z_{k}|,R\}\},
\end{align*}
for every $R>0$. Let $\seq{\mathsf{S}}$ be a sequence of maps of the form
\begin{align}\label{e:weak shifts}
\mathsf{S}_n\co=(z_2,\hdots,z_k,a_n z_1+p_n(z_2,\hdots, z_k))
\end{align}
 such that 
\begin{itemize}
\item[(i)] the degree of each $a_n z_1+p_n(\ze)=d_n \le \tilde{d}$, $\ze \in \C^{k-1}$. In particular, we may assume,
\[p_n(\ze)=\sum_{i \in \mathrm{I}}\alpha_{i,n}\ze^i,\]
where $\mathrm{I}=\{(i_1,i_2,\hdots,i_{k-1}): i_j \in \mathbb{N} \cup \{0\}, 1 \le j \le k-1 \text{ and } 0\le \sum_{j=1}^{k-1}i_j \le \tilde{d}\}$.

\smallskip 
\item[(ii)] there exist uniform constants $\widetilde{m},\widetilde{M}>0$ such that $\widetilde{m}<|a_n|<\widetilde{M}$ and $0\le |\alpha_{i,n}|<\widetilde{M}$ for every $n \ge 1$ and $i \in \mathrm{I}.$
\end{itemize}
The sequence $\seq{\s}$ as above will be called a sequence of \emph{uniformly bounded weak shift-like maps of degree at most $\tilde{d}$}. Further, note that the corresponding non-autonomous sequence of maps defined by
\[S_n\co=\mathsf{S}_n\co+(0,\hdots,0,Q_{d-1}(z_2),\ho_d(z_2,\hdots,z_k))\] 
are of the form (\ref{e:tampered weak-shifts}), if $d \ge \tilde{d}+2$. This shows that the $S_n$'s can be regarded as perturbations of $\mathsf S_n$ provided that $d \geq \tilde d + 2$. Therefore, the sequence $\seq{S}$ will be called a \emph{perturbed sequence of weak shift-like maps of degree $d$} corresponding to the sequence $\seq{\mathsf{S}}$ of weak shift-like maps. 

\smallskip The motivation to consider the above perturbed class of maps is inspired from the study of dynamics of certain (irregular) polynomial automorphisms of $\C^3$ studied in \cite{CF:three dimension}. 

 \begin{lem}\label{l:filtration radius}
Let $\seq{S}$ be a sequence of perturbed weak shift-like maps of degree $d$ corresponding to a uniformly bounded sequence of weak shift-like maps of degree at most $\tilde{d} \le d-2$. Then there exists a sufficiently large $R>1$ such that
\[U^+=\bigcup_{n \ge 1} S(n)^{-1}(\I{V_R^k}) \text{ and }U^{-}=\bigcup_{n \ge 1}\big(S^{-1}(n)\big)^{-1}(\I{V_{R}^-}).\]
\end{lem}
\begin{proof}
Since $\overline{V_R^k}$ contains a neighbourhood of $X^+$ in $\mathbb{P}^k$, the claim would follow if we prove that for every $\ze \in V_R^k$, $S(n,m)(\ze) \to X^+$ as $n \to \infty$ where 
\[
S(n,m):=S(n+m)S(m)^{-1}=S_{n+m}\circ S_{n+m-1}\circ \cdots \circ S_{m+1}.
\] 
For $\ze=\co$, let $\ze'=(z_2,\hdots,z_k)$ and $m_0=\sharp{\mathrm{I}}$, where $\mathrm{I}$ is the set of indices as observed in the definition of weak shift-like maps in (\ref{e:weak shifts}). For $1 \le i \le k$, let $\pi_i$ denote the projection to the $i$-coordinate. Then for every $n \ge 1$ 
\[
|a_n z_1+p_n(\ze')| \le \widetilde{M}\big(m_0+1\big)|z_k|^{d-2}.	
\]
Thus, in $V_R^k$, the growth of the term $a_n z_1+p_n(\ze')$ is bounded by that of $z_k^{d-2}$. Also, we observe below that the growth of the homogeneous term $\ho_d(\ze')$ is bounded by a constant times the growth of $z_k^d$ in $V_R^k$, namely,
\begin{align*}
\frac{|z_k|^d}{k-1}=\Big(1-\frac{k-2}{k-1}\Big)|z_k|^d < |z_k|^d-\frac{k-2}{(k-1)^d}|z_k|^d<|\ho_d(\ze')|<(k-1)|z_k|^d.
\end{align*}
Hence, for $|z_k|$ sufficiently large, i.e., for $R$ sufficiently large
\[|\ho_d(\ze')|>|a_nz_1+p_n(\ze')|\]
whenever $\ze \in V_R^k$. In particular, there exist constants $0<\mathsf{m}<1< \mathsf{M}$ such that 
\begin{align}\label{e:filtration}
\me |z_k|^d< |\pi_k(S_n(\ze))|<\Me |z_k|^d	
\end{align}
for every $n \ge 1$. Next, we will estimate the order of growth of the coordinates $\pi_i \circ S_n(\ze)$, $1 \le i \le k-1$ in $V_R^k$. So, we begin with a further modification of the choice of $R>1$ such that 
\[|\pi_k ( S_n(\ze))| \ge (\me|z_k|)|z_k|^2 \ge (k-1)R^2\] and
\begin{align}\label{e:k-1 term}
\nonumber \me |z_k|^d&\ge \me |z_k||z_k|^{d-1}>\frac{\me}{2} |z_k|\big(|z_k|+|Q_{d-1}(z_2)|\big) \\&
\ge \frac{\me}{2} |z_k||\pi_{k-1}(S_n(\ze))| \ge (k-1)|\pi_{k-1}(S_n(\ze))|.
\end{align} 
Thus, from above and (\ref{e:filtration}),
\[
|\pi_k(S_n(\ze))| > (k-1)|\pi_{k-1}(S_n(\ze))|.
\] 
Further, as $\pi_i \circ S_n(\ze)=z_{i+1}$ for every $1 \le i \le k-2$ it follows that
\begin{align}\label{e:i terms}
|\pi_k(S_n(\ze))|\ge (\me|z_k|)|z_k|^{d-1}>(k-1)^2|\pi_i(S_n(\ze))|
\end{align}
for every $n \ge 1$. Hence, $S_n(V_{R}^k) \subset \I{V_{R^2}^k}.$ Now by induction, it follows that for every $n,m \ge 1$ 
\begin{align}\label{e:growth}
S(n,m)(V_R^k)=S_{m+n}\circ \cdots\circ S_{m+1}(V_R^k) \subset V_{R^{2^n}}^k\subset \I{V_R^k}.	
\end{align}
This shows that for $\ze \in V_{R}^k$, $S(n,m)(\ze) \to \infty$ as $n \to \infty$. 

\medskip \noindent Next, we use the above bounds to observe the dynamics in homogeneous coordinates. Let 
\[{\bf w}_n=[w_1^n:w_2^n:\cdots:1:w_{k+1}^n]=\overline{S(n,m)}([\ze:1]).\]
As a consequence of  (\ref{e:k-1 term}) and (\ref{e:i terms}), 
\[
0 \le |w_i^n| \le \frac{2}{\me R^{2^{n-1}}}
\]
for every $1 \le i \le k-1$. Also, $0<|w_{k+1}^n|<R^{-2^n}.$ Hence ${\bf w}_n \to X^+$ as $n \to \infty$, which completes the proof for $U^+$.

\medskip 

\noindent For the sequence $\{S^{-1}_n\}$, recall that 
{\small \[
S_n^{-1}\co=\Big(a_n^{-1}\big(z_k+(p_n+\ho_d)(z_1,\hdots, z_{k-1}-Q_{d-1}(z_1)\big),z_1,\hdots,z_{k},z_{k-1}-Q_{d-1}(z_1)\Big).\]}
Thus, there exists $\w{\me}<1<\w{\Me}$, such that 
\begin{align}\label{e:negative growth}
    \w{\me} |z_1|^{d^2-d} \le |\pi_1 \circ S_n^{-1}(\ze)|\le \w{\Me} |z_1|^{d^2-d} \text{ and }\w{\me} |z_1|^{d-1} \le |\pi_k \circ S_n^{-1}(\ze)|\le \w{\Me} |z_1|^{d-1}
\end{align}
for $R>0$ sufficiently large. Since $d^2-d \ge d-1$ for $d \ge 3$, 
\[
\norm{S_n^{-1}(\ze))}_{\infty}=|\pi_1 \circ S_n^{-1}(\ze)|
\] 
for every $n \ge 1$ and $R$ large enough. Here, $\supnorm{\cdot}$ denotes the supremum norm on $ \C^k$. Thus, $S_n^{-1}(V_R^-) \subset V_{R^2}^-$ and the order of growth of the first coordinate is larger than the other coordinates. Hence, by similar arguments as in the proof of $U^+$, it follows that for every $m\ge 1$, $S^{-1}(n,m)(\ze) \to X^-$, whenever $n \to \infty$ and $\ze \in U^-$. 
\end{proof}

\begin{rem}\label{r:nested}
Note that (\ref{e:growth}) also shows that $S(n)^{-1}(V_R^k) \subset S(n+1)^{-1}(V_R^k)$ for every $n \ge 1.$ and hence $U^+$ is an open connected subset of $\C^k.$
\end{rem}

\begin{rem}\label{r:all filtration}
For $2 \le  i\le  k$, consider the sets 
\[V_{R}^i=\big\{z \in \C^k:|z_i| \ge \max\{(k-1)|z_j|,R: 1\le j \le k, i \neq j\}\big \} \text{ and } V_R^+=\bigcup_{i=2}^{k}V_{R}^i.\]
Then, for $z \in V_R^i$ the bounds in (\ref{e:filtration}) and (\ref{e:k-1 term}) in the proof of Lemma \ref{l:filtration radius} change as follows:
\begin{align*}
\me |z_i|^d< &|\pi_k(S_n(\ze))|<\Me |z_i|^d	\\
\me |z_i|^d\ge \me |z_i||z_i|^{d-1}>\frac{\me |z_k|}{2}\big(&|z_k|+|Q_{d-1}(z_2)|\big)\ge (k-1)|\pi_{k-1}(S_n(\ze))|.		
\end{align*}
Thus, $S_n(V_R^+) \subset V_{R^2}^k$ and (\ref{e:filtration}) can be stated as -- for every $z \in V_R^+$ and $n \ge 1$
\begin{align}\label{e:final filtration}
\me\supnorm{\ze}^d \le \supnorm{S_n(\ze)}=\abs{\pi_k\circ S_n(\ze)} \le \Me \supnorm{\ze}^d
\end{align}
where, as before, $\supnorm{\cdot}$ denotes the supremum norm on $\C^k$.
\end{rem}

\begin{rem}\label{r:modified radius}
Let \[W_R^-:=\{z \in \C^k:(k-1)|z_1| \ge  \max\{|z_2|,\hdots, |z_{k}|,R\}\}.\]
Then, as a consequence of (\ref{e:negative growth}), the filtration radius $R>0$ which has been chosen in Lemma \ref{l:filtration radius}, can be further appropriately increased so that $S_n^{-1}(W_R^-)=V_{R^2}^-$ for every $n \ge 1$. Though the set $W_R^-$ may not be relevant immediately, we will need this observation later in Section \ref{s:4}, particularly in the proof of Lemma \ref{l:phi onto}. 
\end{rem}
Consider the sequence 
\[G_n(z)=\frac{1}{d^n} \log^+\supnorm{S(n)(z)}\]
of plurisubharmonic functions on $\C^k$. Note that $G_n$ is pluriharmonic on $V_R^+$ for every $n \ge 1$. The goal now will be to show that the $G_n$'s converge uniformly on compact subsets of $\C^k$ to a non-constant continuous plurisubharmonic function $G$, we call it {\it the non-autonomous Green's function corresponding to the sequence $\seq{S}$}. This will be done by first showing that the limit exists on compact subsets of $U^+$ and then extending the domain of convergence to all of $\C^k$. All this forms the content of Proposition $3.7$, Theorem $3.8$ and Lemma $3.9$.

\medskip

In Section \ref{s:4}, it will be important to know that the same limit $G$ is achieved when working in the non-autonomous setting. To this end, the following definitions will be required. Corresponding to a given sequence $\seq{S}$ and every fixed $m \ge 1$ we define a new sequence as follows
\begin{align*}
\mathbf{S}^m_n=	S_{n \text{ mod } m }\text{ and } s_m=S(m)=S_m \circ \cdots \circ S_1
\end{align*}
with the convention that $lm \text{ mod }m=m$ for every integer $l\ge 1$, i.e.,
\[\mathbf{S}^m_1=S_1, \mathbf{S}^m_2=S_2,\hdots, \mathbf{S}^m_m=S_m, \mathbf{S}^m_{m+1}=S_1,\mathbf{S}^m_{m+2}=S_2, \hdots , \mathbf{S}^m_{2m}=S_m,\hdots  \]
Thus, for $l \ge 1$,
\[
S(lm)=\overbrace{S(m)\circ \cdots \circ S(m)}^{l-\text{times}}=s^l_m.
\]
Note that for every $m \ge 1$, $\seq{\mathbf{S}^{m}}$ is a non-autonomous sequence satisfying the assumptions (i) and (ii) stated just after $(3.2)$. Let $\G_m$, $U_m^+$ and $K_m^+$ denote the non-autonomous Green's function, the escaping set and the non-escaping set respectively, corresponding to the sequence $\seq{\mathbf{S}^{m}}.$ Also, let $\mathcal{K}_m^+=\C^k \setminus U_m^+.$ Now note that from the dynamical point of view, the non-autonomous dynamics of the sequence $\seq{\mathbf{S}^m}$ is equivalent to the iterative dynamics of the map $s_m$. In particular, note that 
\[\G_m(z)=\lim_{n \to \infty}\frac{1}{d^n}\log^+\supnorm{\mathbf{S}^m_n(z)}=\lim_{l \to \infty}\frac{1}{d^{lm}}\log^+\supnorm{S(lm)(z)}=\lim_{l \to \infty}\frac{1}{d^{lm}}\log^+\supnorm{s_m^l(z)}.\]
Hence
\begin{align}\label{e:iterative Green}
\G_m\big(s_m(z)\big)=d^m \G_m(z)
\end{align}
for every $m \ge 1$. In Lemma $3.10$, it will be shown that the $\G_m$'s also converge to $G$ uniformly on compact subsets of $\C^k$. But we first start by studying the $G_n$'s.

\begin{prop}\label{p:Green sequence}
The sequence $\seq{G}$ converges on compact subsets of $U^+$ to $\w{G}$, a continuous, non-constant pluriharmonic function on $U^+$.
\end{prop}
\begin{proof}
As in the case of H\'{e}non maps, we will use the filtration bounds (\ref{e:filtration}). Let $C$ be a compact subset of $U^+$ and let $V^+=V_{R}^+$, where $R>1$ is appropriately chosen as observed in Lemma \ref{l:filtration radius}. 

\medskip

\textit{Step 1:} There exists a positive integer $N_{C} \ge 1$ such that $S(n)(C) \subset {V^+}$ for every $n \ge N_C.$

\smallskip\noindent For $z \in C$, Lemma \ref{l:filtration radius} shows that there exists $n_z \ge 1$ such that $S(n_z)(z) \in \I{V_R^k} \subset V^+$. Also, by Remark \ref{r:nested}, $S(n)(z)\in V^+$ for every $n \ge n_z.$ Since $S(n)$ is a biholomorphism for every $n \ge 1$, there exists a ball of radius $\delta_z>0$ such that $S(n)(B(z;\delta_z)) \subset V^+$ for every $n \ge n_z.$ As $C$ is compact in $U^+$,  there exists a finite collection of balls $B(z_i;\delta_{z_i})$, $1 \le i \le n_C$ that cover $C$. There exist corresponding $n_{z_i}$'s such that $S(n)(B(z_i,\delta_{z_i})) \subset V^+$ whenever $n \ge n_{z_i}.$ Thus for $N_C=\max\{n_{z_i}: 1 \le i \le n_C\}$, the claim in \textit{Step 1} holds. 

\medskip

\textit{Step 2:} For every $w \in V^+$, there exists a neighbourhood of $w$, say $B(w;\delta_w)$, $\delta_w>0$, on which the sequence $\seq{G}$ is uniformly Cauchy.

\smallskip\noindent Note that $\log^+\supnorm{z}=\log\supnorm{z}>0$ for every $z \in B(w:\delta_w) \subset V^+.$  By (\ref{e:filtration}), (\ref{e:growth}) and (\ref{e:final filtration}),
\begin{align}\label{e:recursive}
\frac{1}{d^{n+1}}\log \me +\frac{1}{d^n}\log^+\supnorm{S(n)(z)}&\le \frac{1}{d^{n+1}}\log^+\supnorm{S(n+1)(z)} \\
&\le \frac{1}{d^{n+1}}\log \Me + \frac{1}{d^n}\log^+\supnorm{S(n)(z)},	
\end{align}
for every $n \ge 0$, where $0<\me<1<\Me$ are as in Lemma \ref{l:filtration radius}. Hence
\begin{align}\label{e:non-constant}
-\sum_{i=1}^n\frac{\widetilde{\Me}}{d^i}+\log^+\supnorm{z} \le G_n(z) \le \sum_{i=1}^n\frac{\widetilde{\Me}}{d^i}+\log^+\supnorm{z},	
\end{align}
 where 
$\widetilde{\Me}=\max\{|\log \me|,|\log \Me|\}$. Thus, $\seq{G}$ is uniformly bounded in $B(w;\delta_w)$. Further, from (\ref{e:recursive})
\[|G_n(z)-G_{n+1}(z)| \le \frac{\widetilde{\Me}}{d^{n+1}}\] and this shows that $\seq{G}$ converges uniformly to a function $\w{G}$ on $B(w;\delta_w)$, thus completing \textit{Step 2}.

\medskip 
Since $G_n(z) - \log^+\supnorm{z}$ is a bounded pluriharmonic function on $V^+$ for every $n \ge 1$, we have $G_n(z) = \log^+\supnorm{z}+C_n$, where $\seq{C}$ is a convergent sequence. Thus $\seq{G}$ converges locally uniformly to a continuous, non-constant, pluriharmonic function $\w{G}$ on $V^+$. 

\medskip
\textit{Step 3:} The sequence $\seq{G}$ converges uniformly on compact subsets of $U^+$.

\medskip Given a compact set $C \subset U^+$ consider the non-autonomous sequence $\seq{S^C}$ such that $S^C_n=S_{N_C+n}$, where $N_C\ge 1$ is as obtained in \textit{Step 1}. Define the sequence of plurisubharmonic functions on $\C^k$ as
\[
G^C_n(z)=\frac{1}{d^n}\log^+\supnorm{S^C(n)(z)}=\frac{1}{d^n}\log^+\supnorm{S(n, N_C)(z)}
\] 
which are pluriharmonic on $V^+$. Thus by \textit{Step 2}, $\{G^C_n\}$ converges uniformly on all compact subsets of $V^+$ to a function $G^C$. Note that by definition, for every $z \in \C^k$ and $n \ge 1$
\begin{align}\label{e:convergence on compact}
G^C_n\circ S(N_C)(z)=d^{N_C}G_n(z), \text{ i.e., }G^C_n(z)=d^{N_C}G_n\circ S(N_C)^{-1}(z).
\end{align}
Now, from \emph{Step 2}, both the sequences $\{G^C_n\}$ and $\{G_n\}$ converge uniformly on compact subsets of $V^+$. Since $S(N_C)(C)$ is a compact set contained in $V^+$, from (\ref{e:convergence on compact}) above, it follows that for $\{G_n\}$ is uniformly convergent on $C$, which proves \emph{Step 3}.

\medskip Further, note that $G_n^C \circ S(N_C)$ is pluriharmonic on a neighbourhood of $C$, in fact on the finite cover of $C$ as observed in \emph{Step 1} and hence $\tilde{G}$ is pluriharmonic on $U^+$.
\end{proof}

Next, extend $\w{G}$ to $\C^k$ as
\begin{align}\label{e:Green}
G(z)=\begin{cases}
	\w{G}(z) &\text{ for } z \in U^+, \\
	0 &\text{ for } z \in \mathcal{K}^+=\C^k \setminus U^+
\end{cases}
\end{align}
and with this, we generalise Theorem 2.2 from \cite{SG:regular} as 
\begin{thm}\label{t:Green function}
The function $G$ in (\ref{e:Green}) is continuous, non-constant and plurisubharmonic on $\C^k$.
\end{thm}
\begin{proof}
Note that to prove continuity of $G$ on $\C^k$, we only need to prove its continuity at points on $\partial \mathcal{K}^+$. Suppose not and assume that $G$ is discontinuous at $z_0 \in \partial \mathcal{K}^+$. In particular, there exists a  sequence $\seq{z} \in U^+$ such that $z_n \to z_0$ but $G(z_n)>c$ for some $c>0$. 

\medskip 
Our goal will be to contradict the above assumption. The idea behind it is to show that any point $z \in \C^k$ with $G(z)>c>0$ should diverge to infinity under the non-autonomous action of the family $\{S_n\}$ approximately at an order $c^{d^n}$. However, as a consequence of Lemma \ref{l:filtration radius} and its proof, it turns out that the a convergent sequence of points with the former property should be contained in $U^+$, i.e., eventually contained in $V_{\w{R}}^+$ uniformly, where $\w{R}>0$ is sufficiently large and appropriately chosen. This will force $z_0 \in \mathcal{K}^+$, as chosen in the above paragraph, to be contained in $U^+$, thus giving a contradiction!

\medskip To begin with, recall the constant $\w{M}$ and the index set $\mathrm{I}$, from the assumptions (i) and (ii) on the sequence $\seq{S}$. Since $z_n \to z_0$, $\supnorm{z_0} \le r_{\seq{z}}$ where $r_{\seq{z}}=\sup\{\supnorm{z_n}: n \ge 1\}<\infty$. Let
\[\w{R}>\max \{r_{\seq{z}}, 2(k-1), (m_0+k)R^d\}\] 
where $R>1$ as in Lemma \ref{l:filtration radius} and $m_0=\sharp \mathrm{I}$. Further, redefine the constant $\w{\Me}>0$ as
\[\w{\Me}=\max \{\log \w{R}, |\log \Me|, |\log \me|, |\log\big(\w{M}(m_0+k)\big)|\}\]
and let $\ku_0 \ge 1$ such that for $m \ge \ku_0$ 
\[\frac{\w{\Me}}{d^m}<\frac{c}{8} < c \text{ and }\Big(\frac{d-1}{d}\Big)^{m-1}\log\w{R}<\frac{c}{8}.\] 

\medskip
\noindent \textit{Claim}: $S(m)(z_n) \in \I{V_{\w{R}}^+}$ for every $m \ge \ku_0$ for all $n \ge 1$.

\medskip 

 Suppose not. Then there exists  a subsequence $\{z_{n_l}\} \subset \seq{z}$ and a sequence of natural numbers $\{m_l\}$ with each $m_l \ge \ku_0$ and such that $S(m_l)(z_{n_l}) \notin \I{V_{\w{R}}^+}$. By abuse of notation, let $z_l$ be any point in the sequence $\{z_{n_l}\}$ and let $m_l \ge \ku_0$ be the corresponding positive integer such that $S(m_l)(z_{l}) \notin \I{V_{\w{R}}^+}$.

 \smallskip
 Also, from (\ref{e:growth}) and (\ref{e:final filtration}) it follows that if 
 $S(\widetilde{m}_l)(z_l) \in {V_{\w{R}}^+}$ then  $S(m)(z_l) \in V_{\w{R}}^+$ for every $m \ge \widetilde{m}_l$. Hence it follows that
 \[ S(m)(z_l) \notin {V_{\w{R}}^+} \text{ whenever }m <  m_l.\]
Now since $z_l \in U^+$, there exists  a minimum integer $\ku_l\ge m_l \ge \ku_0$ such that 
\[S(m)(z_l) \in V_{\w{R}}^+ \text{ for } m \ge \ku_l \text{ and } S(m)(z_l) \notin V_{\w{R}}^+ \text{ for } m < \ku_l.\]
\noindent Further, if $S(\ku_l-1)(z_l) \in V_{\w{R}}^-$, then by using Lemma \ref{l:filtration radius}, it follows that
 \[z_l \in S(\ku_l-1)^{-1}(V_{\w{R}}^-) \subset V_{\w{R}^{2^{\ku_l-1}}}^-,\] 
 i.e., $\supnorm{z_l}>\w{R}^{2}>\w{R}$ which is a contradiction! In particular, we can now assume that 
 \[S(m)(z_l) \notin \overline{ V_{\w{R}}^+ \cup V_{\w{R}}^-}, \text{ i.e., } S(m)(z_l) \notin {\bigcup_{i=1}^k V_{\w{R}}^i}\] 
 for every $1\le m \le \ku_l-1.$ Also, as before, since $z_l \in U^+$ there exists $p_l \ge 0$, the minimum possible value, such that 
 \[\supnorm{S(m)(z_l)}>1 \text{ for every }m\ge p_l \text{ and } \supnorm{S(p_l-1)(z_l)} \le 1<R<\w{R}.\] 
 Let $w_0=S(p_l)(z_l)$ and $w_0=z_l$ if $p_l=0$. Note that by the choice of $\w{R}>0$ it follows that $0 \le p_l\le \ku_l-1$ and $\supnorm{w_0}\le \w{R}.$ 
 Let 
 \[w_{m}=S_{p_l+m}(w_{m-1})\text{ whenever }1\le m\le \ku_{l}-p_l-1.\]
 By the assumptions, $w_m \notin V_{\w{R}}^i$ for every $1 \le i \le k$ and $1 \le m \le \ku_l-p_l-1$. Also $\supnorm{w_m} >1$. Now, if $|\pi_i(w_{m+1})| =\supnorm{w_{m+1}}$ for some $1 \le i \le k-2$, then $\supnorm{w_{m+1}}<\supnorm{w_m}$. If $\supnorm{w_{m+1}}=|\pi_{k-1}(w_{m+1})|$ then
 \[\supnorm{w_{m+1}}=|\pi_{k-1}(w_{m+1})| \le  2\supnorm{w_m}^{d-1}\le 2(k-1)\supnorm{w_m}^{d-1}.\] 
 Otherwise if $\supnorm{w_{m+1}}=|\pi_{k}(w_{m+1})|$ then either $\supnorm{w_m}\le \w{R}$ or there exists $1 \le i \le k-1$ such that
 \[\supnorm{w_{m+1}}=|\pi_k(w_{m+1})| \le (k-1)|\pi_{i}(w_{m+1})|\le 2(k-1)\supnorm{w_m}^{d-1}.\] 
 So, for $1 \le m \le \ku_l-p_l-1$,
 \[\supnorm{w_{m}} \le \big(2(k-1)\big)^{\sum_{i=0}^{m-1} (d-1)^i}\w{R}^{(d-1)^m},\]
 i.e., 
 \begin{align*}
 \log^+\supnorm{w_{m}} &\le \log \big(2(k-1)\big)\sum_{i=0}^{m-1} (d-1)^i+{(d-1)^m}\log \w{R},
 \\ &\le\log \big(2(k-1)\big)\sum_{i=0}^{m+p_l-1} (d-1)^i+{(d-1)^{m+p_l}}\log \w{R}
 \end{align*}
 Hence for $m=\ku_l-p_l-1$ we have $w_{\ku_l-p_l-1}=S(\ku_l-1)(z_l)$ and
 \begin{align*}
 \log^+\supnorm{S(\ku_l-1)(z_l)} \le \log\big(2(k-1)\big)\sum_{i=0}^{\ku_l-2} (d-1)^i+{(d-1)^{\ku_l-1}}\log \w{R}.
 \end{align*}
 Thus
 \begin{align*}
\frac{1}{d^{\ku_l-1}}\log^+\supnorm {S(\ku_l-1)(z_l)} &\le \Big(\frac{d-1}{d}\Big)^{\ku_l-1}\Bigg(\log\w{R}+\log 2(k-1) \sum_{i=1}^{\ku_l-1} \frac{1}{(d-1)^i}\Bigg) \\	
&\le \Big(\frac{d-1}{d}\Big)^{\ku_l-1}2\log\w{R},
 \end{align*}
and 
\begin{align*} 
\frac{1}{d^{\ku_l}}\log^+\supnorm {S(\ku_l)(z_l)}  &\le\frac{1}{d^{\ku_l}}\log \w{M}(m_0+k) + \frac{1}{d^{\ku_l-1}}\log^+\supnorm{S(\ku_l-1)(z_l)} \\& \le \Big(\frac{d-1}{d}\Big)^{\ku_l-1}2\log\w{R}+\frac{\w{\Me}}{d^{\ku_l}.}
\end{align*}

 Also, $S(\ku_l)(z_l) \in V_{\w{R}}^+$ and satisfies (\ref{e:final filtration}). By the choice of $\w{\Me}$,
 \begin{align*} 
 \frac{1}{d^{\ku_l+1}}\log^+\supnorm{S(\ku_l+1)(z_l)} &\le \frac{1}{d^{\ku_l}}\log^+\supnorm {S(\ku_l)(z_l)}+\frac{\w{\Me}}{d^{\ku_l+1}}\\& \le  \Big(\frac{d-1}{d}\Big)^{\ku_l-1}2\log\w{R}+\frac{\w{\Me}}{d^{\ku_l}}+\frac{\w{\Me}}{d^{\ku_l+1}}.
 \end{align*}
 Thus, for $m \ge \ku_l+1$
 \begin{equation*}
0\le \frac{1}{d^m}\log^+\supnorm{S(m)(z_l)} \le \Big(\frac{d-1}{d}\Big)^{\ku_l-1}2\log\w{R}+\sum_{i=\ku_l}^m \frac{\w{\Me}}{d^i}
< \frac{c}{4}+\frac{c}{4}<\frac{c}{2}
 \end{equation*} 
 by the choice of $\ku_l\ge \ku_0$. Hence, $G(z_l)\le c/2$, and this contradicts the assumption on the sequence $\seq{z}$. This also proves the claim.
 
 \medskip Now, by Remark \ref{r:all filtration}, $S(\ku_0+1)(z_n) \in V_{\w{R}}^k$ for every $n \ge 1.$ Now since $z_n \to z_0$ and $\w{R}>R$, by continuity $S(\ku_0+1)(z_0) \in V_{\w{R}}^k \subset {V_{R}^k}$. Hence, from the proof of Lemma \ref{l:filtration radius}, $z_0 \in U^+$ and $z_0 \notin \partial \mathcal{K}^+ \subset \mathcal{K}^+$, which is a contradiction to the assumption.

 \medskip  Note that the function $G$ as defined in (\ref{e:Green}) is pluriharmonic on both $\I{\mathcal{K}^+}$ and $U^+$. Since $\C^k=\mathcal{K}^+ \cup U^+$, the upper semicontinuous regularization of $G$ (on $\partial \mathcal{K}^+$ from $U^+$ and $\I{\mathcal{K}^+}$) will be a plurisubharmonic on $\C^k$. Now, as $G$ is continuous it coincides with the upper semicontinuous regularization and is plurisubharmonic. 
 \end{proof}
 
 \begin{lem}\label{l:uniform convergence 1}
 The sequence $\seq{G}$ converges uniformly to $G$ on compact subsets of $\C^k$.
 \end{lem}
\begin{proof}
The proof is divided into three cases. Let $C$ be a compact subset of $\C^k$.

\medskip
\noindent \textit{Case 1: }If $C \subset U^+$, the proof is immediate from Proposition \ref{p:Green sequence}. 

\medskip

\noindent \textit{Case 2: }Suppose $C \subset \I{\mathcal{K}^+}$. Then, as in the proof of Theorem \ref{t:Green function}, consider
\[\w{R}>\max\{ 2(k-1), (m_0+k)R^d,\supnorm{z}: z \in C\}\]
where $R$ is the radius of filtration obtained in Lemma \ref{l:filtration radius}. Also, by the same arguments as in the proof of Theorem \ref{t:Green function}, $S(n)(z) \notin V_{\w{R}}^-$ for every $n \ge 1.$ Thus, for $z \in C$, $S(n)(z) \in \overline{(V_{\w{R}}^+ \cup V_{\w{R}}^-)^c}$ or $\norm{S(n)(z)} \le \w{R}$ for every $n \ge 1$ and $z \in C.$ 
Further by similar arguments as in the proof of Theorem \ref{t:Green function}, it follows that the order of growth of the sequence $\{S(n)(z)\}$ is at most exponential of ${(d-1)^n}$ for some sufficiently large enough constant. In particular, if $\supnorm{S(n-1)(z)}>1$ then
\[\supnorm{S(n)(z)} \le  2(k-1)\supnorm{S(n-1)(z)}^{d-1} .\]
Hence inductively for every $z \in C$
\[\supnorm{S(n)(z)}\le \big(2(k-1)\big)^{\sum_{i=1}^{n-1}(d-1)^i}\w{R}^{(d-1)^n}.\]
Thus, for a given $\ep>0$ and for every $z \in C$ there exists a sufficiently large $n_0 \ge 1$ such that
\[0 \le G_n(z) \le \Big(\frac{d-1}{d}\Big)^{n}\Bigg(\log\w{R}+\log 2(k-1) \sum_{i=1}^{n-1} \frac{1}{(d-1)^i}\Bigg) < \Big(\frac{d-1}{d}\Big)^{n}2\log\w{R} \le \ep.\]
whenever $n \ge n_0.$

\medskip 

\noindent \textit{Case 3:} Suppose that both $C_0=C \cap U^+ \neq \emptyset$ and $\w{C}=C \cap \mathcal{K}^+ \neq \emptyset$. In this case, by \textit{Case 2}, for a given $\ep>0$ there exists $\ku_1 \ge 1$ such that $|G_n(z)-G(z)|<\ep$ for $n \ge \ku_1$ and $z \in \w{C}$.  

\smallskip\noindent Define the sequence of subsets of $C_0$ as
\[ C_n=C_0 \cap \Big(S(n)^{-1}\big(\I{V_{\w{R}}^+}\big)\setminus S(n-1)^{-1}\big(\I{V_{\w{R}}^+}\big)\Big).\]
Note that $C_0 \cap V_{\w{R}}^+ \subset C_0 \cap S_1^{-1}(V_{\w{R}}^+)=C_1$. Hence $\cup_{n=1}^\infty C_n=C_0.$ By continuity of $G$, we can find $\ku_2 \ge \ku_1$ such that for every $z \in C_n$ and $n \ge \ku_2$, $|G(z)| \le \frac{\ep}{2}$. Now, let $\ku_3 \ge \ku_2$ be such that for $n \ge \ku_3$
\[ \frac{\w{\Me}}{d^n} < \frac{\ep}{8} \text{ and } \Big(\frac{d-1}{d}\Big)^{n-1}\log\w{R}<\frac{\ep}{8},\]
where $\w{\Me}$ is the constant as constructed in the proof of Theorem \ref{t:Green function}, i.e.,
 \[\w{\Me}=\{\log \w{R}, |\log \Me|, |\log \me|, |\log \w{M}(m_0+k)|\}.\]
Suppose $z \in C_n$, $n \ge \ku_3$ then by the same argument as in $\textit{Case 2}$, it follows that
\[ \frac{1}{d^{n-1}}\log^+\supnorm{S(n-1)(z)}  \le \Big(\frac{d-1}{d}\Big)^{n-1}2\log\w{R} \text{ and } S(n)(z) \in V_{\w{R}}^+.\] 
Further, by the assumption on the sequence $\seq{S}$, 
\begin{align*}
\frac{1}{d^n}\log^+\supnorm{S(n)(z)}  &\le \Big(\frac{d-1}{d}\Big)^{n-1}2\log\w{R} +\frac{1}{d^n}\log\w{M}(m_0+k)\\  &\le  \Big(\frac{d-1}{d}\Big)^{n-1}2\log\w{R}+\frac{\w{\Me}}{d^n}.
\end{align*}
Thus, for $m \ge 1$ 
\begin{align*}
\frac{1}{d^{n+m}}\log^+\supnorm{S(n+m)(z)} &\le \Big(\frac{d-1}{d}\Big)^{n-1}2\log\w{R}+\sum_{i=n}^m \frac{\w{\Me}}{d^i}
< \frac{\ep}{4}+\frac{\ep}{4} < \frac{\ep}{2}
\end{align*}
and this means that for $z \in C_n$, $n \ge \ku_3$
\[|G(z)-G_n(z)| \le |G(z)|+|G_n(z)|<\ep.\]
Let $C'=\cup_{i=1}^{\ku_2-1} C_i$ and $C''=\cup_{i=\ku_2}^\infty C_i$. As $C'$ is compactly contained in $U^+$, by \textit{Case 1}, there exists $\ku_4 >\ku_3$ such that
\[|G_n(z)-G(z)|<\ep \text{ for every }z \in C'.\]
Note that $C=C' \cup C'' \cup \w{C}$, thus $|G_n(z)-G(z)|<\ep$ for every $n \ge \ku_4$ for every $z \in C.$
\end{proof}

Recall the sequence $\{\mathbf S^m_n\}$ and the dynamical objects associated with it, as defined just before Proposition $3.7$. 
\begin{lem}
For $\ep > 0$ and $C \subset \C^k$ a compact subset, there exists $\ku_C \ge 1$ such that for $m \ge \ku_C$,
$|G-\G_m|_C<\ep$.
\end{lem}
 
 \begin{proof}
 Since the idea of the proof is similar to that of Lemma \ref{l:uniform convergence 1}, we will be brief. Let $\w{\Me}$ and $\w{R}$ be as chosen there.
 
\medskip
\noindent \textit{Case 1: }Suppose $C \subset U^+$. Then by Proposition \ref{p:Green sequence} there exists $\ku_0>1$, large enough, such that $S(m)(C) \subset V_{R}^+$ for $m \ge \ku_0.$ Thus $C \subset U_m^+$ for $m \ge \ku_0.$ Further, from (\ref{e:filtration}) and (\ref{e:growth}) it follows that for every $l \ge 2$ 
 \begin{align*} 
 -\sum_{i=m+1}^{lm} \frac{\w{\Me}}{d^i}+\frac{1}{d^m}\log^+\supnorm{s_{m}(z)} &\le \frac{1}{d^{lm}}\log^+\supnorm{s_{m}^l(z)}\le \frac{1}{d^{m}}\log^+\supnorm{s_{m}(z)} +\sum_{i=m+1}^{lm} \frac{\w{\Me}}{d^i},
 \end{align*}
i.e.,
 \[ -\sum_{i=m+1}^{lm} \frac{\w{\Me}}{d^i}+G_m(z)\le \frac{1}{d^{lm}}\log^+\supnorm{s_{m}^l(z)} \le G_m(z)+\sum_{i=m+1}^{lm} \frac{\w{\Me}}{d^i}.\]
 Hence, for $m$ sufficiently large, $|G_m-\G_m|<\w{\Me}/{d^m}< \ep $ on $C$. Thus, by Lemma \ref{l:uniform convergence 1}, the result holds in this case.
 
\medskip

\noindent \textit{Case 2: }Suppose $C \subset \mathcal{K}^+$. Let $\ep>0$. Choose $\ku_0 \ge 1$ such that for $m \ge \ku_0$
 \[ \frac{\w{\Me}}{d^m} < \frac{\ep}{8} \text{ and } \Big(\frac{d-1}{d}\Big)^{m-1}\log\w{R}<\frac{\ep}{8}.\]
 Then $s_m(z)=S(m)(z) \notin \overline{V_{\w{R}}^+ \cup V_{\w{R}}^-}$ for every $m \ge 1$ and $z \in C.$ Thus,
 \[ \frac{1}{d^{m-1}} \log^+\supnorm{s_m(z)}  \le \Big(\frac{d-1}{d}\Big)^{m-1}2\log\w{R}.\]
 
\medskip
 
\noindent \textit{Subcase 1: }Suppose there exists $z \in C$ such that $z \in U^+_m$ for some $m \ge \ku_0$.

\smallskip Choose the smallest $l \ge 1$ such that $s_m^{l+1}(z) \in V_{\w{R}}^+.$ By definition, there is a largest value $l_0$ with $1 \le l_0 < m$ such that 
\[S_{l_0} \circ S_{l_0-1}\circ \cdots S_1\circ s_m^l(z) \notin V_{\w{R}}^+.\] 
Thus, by the assumption on the sequence $\seq{S}$ and arguments as before
\[\frac{1}{d^{lm+l_0+1}} {\log^+\supnorm{S_{l_0+1} \circ S_{l_0}\circ \cdots S_1\circ s_m^l(z)}} \le \Big(\frac{d-1}{d}\Big)^{lm+l_0}\log\w{R}+\frac{\w{\Me}}{d^{lm+l_0+1}}.\]
As $lm+l_0 \ge m \ge \ku_0$,
\begin{align*}
	0\le \frac{1}{d^{nm}} \log \supnorm{s_m^n(z)} &\le \Big(\frac{d-1}{d}\Big)^{lm+l_0}\log\w{R}+\sum_{i=lm+l_0+1}^{nm} \frac{\w{\Me}}{d^{i}} <\frac{\ep}{2}.
\end{align*}
Hence for $m \ge \ku_0$, $|\G_m(z)|<\epsilon/2$, in this case.

\medskip
\noindent \textit{Subcase 2: }Suppose there exists $z \in C$ such that $z \in \mathcal{K}^+_m$ for some $m \ge \ku_0$.

\smallskip\noindent Note that then $0=|\G_m(z)|<\ep/2$. Thus $|G_m|<\ep/2$ for every $z \in C \subset \mathcal{K}^+$ and $m \ge \ku_0$, which completes the argument for \textit{Case 2}.

\medskip

\noindent \textit{Case 3: }Suppose that both $C_0=C \cap U^+ \neq \emptyset$ and $\w{C}=C \cap \mathcal{K}^+ \neq \emptyset$. By \textit{Case 2}, for a given $\ep>0$ there exists $ \ku_0 \ge 1$ such that $|\G_m(z)-G(z)|<\ep$ for $m \ge \ku_0$ and $z \in \w{C}$. For $\ku_1> \ku_0$, let
\[ C_{\ku_1}=C_0 \cap \overline{S(\ku_1)^{-1}(V_{\w{R}}^+)} \text{ and } C'_{\ku_1}=C_0 \setminus C_{\ku_1}.\]
Now, we may assume $\ku_1>\ku_0 \ge 1$ is large enough so that $|G(z)|<\ep/2$ on $C'_{\ku_1}$. Since $C_{\ku_1}$ is compact in $U^+$ by \textit{Case 1}, there exists $\ku_2 \ge \ku_1 \ge 1$ such that $|\G_m^+-G|_{C_{\ku_1}}<\ep$ for $m \ge \ku_2.$ Let $z \in C'_{\ku_1}$ and $m \ge \ku_2$.

\begin{itemize}
 \item If $s_m(z) \notin V_{\w{R}}^+$, then by the same argument as in \textit{Case 2}, we have $|\G_m(z)|<\ep/2.$
 \item If $s_m(z) \in V_{\w{R}}^+$ then there exists a maximum $\ku_1<l_0<m$ such that 
 \[ S(l_0)(z) \notin V_{\w{R}}^+ \text{ and }S(l_0+1)(z)\in V_{\w{R}}^+.\]
\end{itemize}
Thus, by the arguments in \textit{Case 2}
\begin{align*}
	0\le \frac{1}{d^{nm}}\log \supnorm{s_m^n(z)} &\le \Big(\frac{d-1}{d}\Big)^{l_0}\log\w{R}+\sum_{i=l_0+1}^{nm} \frac{\w{\Me}}{d^{i}}<
	\frac{\ep}{2}.
\end{align*}
Hence $|\G_m(z)|<\ep/2$ for every $z \in C_{\ku_1}'$, and therefore $|\G_m-G|_C<\ep$.
  \end{proof}
  
 \section{Non-autonomous attracting basins of weak shift-like maps}\label{s:4}
 
\begin{thm}\label{t:na weak-shift basin}
 Let $\seq{\mathsf{S}}$ be a sequence of weak shift-like maps of degree at most $\tilde{d}\ge 1$, satisfying assumptions (i) and (ii). Further, assume that the origin is a uniformly attracting fixed point for every $\mathsf{S}_n$, i.e., there exist constants $A, B$ such that $0<\A<\B<1$ and $\mathsf{r}>0$ such that
 \begin{align}\label{e:attracting-shift}
     \A\supnorm{z}\le |\pi_k \circ \mathsf{S}_n(z)| \le \B\supnorm{z} 
 \end{align}
  for every $n \ge 1$ and $z \in B(0;\mathsf{r})$. Then the basin of attraction of $\seq{\mathsf{S}}$ at the origin
 \[ 
 \Omega_{\seq{\mathsf{S}}} =\{z \in \C^k:\mathsf{S}(n)(z) \to 0 \text{ as }n \to \infty\}
 \]
 is biholomorphic to $\C^k$.
 \end{thm}
  \begin{proof} Choose $d\ge \tilde{d}+3$ large enough so that $\B^{d-2}< {\A}/{k}$. Define the following two sequences of automorphisms of $\C^k$
\[H_n=\mathsf{S}(nk,(n-1)k)=\s_{nk} \circ \s_{nk-1}\circ \cdots \circ \s_{(n-1)k+1}\]and
\[F_n=S(nk,(n-1)k)=S_{nk} \circ S_{nk-1}\circ \cdots \circ S_{(n-1)k+1}\]
where $S_n$ is the degree $d$-perturbation of the sequence $\seq{\s}$, i.e., 
\[S_n\co=\s_n\co+(0,\hdots,Q_{d-1}(z_2),\ho_d(z_2,\hdots,z_k)).\]
Then $F_n(0)=0=H_n(0)$  and $DF_n(0)=DH_n(0)$ for every $n \ge 1$. Also, note
\[\A\supnorm{z}\le \supnorm{H_n(z)} \le \B\supnorm{z}\]
for every $z \in B(0;\mathsf{r})$. In particular, in the Euclidean norm
\[\frac{\A^l}{\sqrt{k}}\norm{z}\le \norm{H\big((n+1)l, nl\big)(z)} \le \sqrt{k}\B^l\norm{z}\]
for every $n,l \ge 1$. Choose $l_0\ge d-2$ such that $\sqrt{k}\B^{l_0}<1.$ Then 
\[
(\sqrt{k}\B^{l_0})^{d-2} \le \frac{\A^{l_0}}{k^{(d-2)/2} k^{l_0-d+2}} \le \frac{\A^{l_0}}{\sqrt{k}}. 
\]
Let $\w{F}_n=F\big((n+1)l_0,nl_0+1\big)$ and $\w{H}_n=H\big((n+1)l_0,nl_0+1\big)$. Then the Taylor expansions of $\w{F}_n$ and $\w{H}_n$ at the origin agree up to order $d-2$ for every $n \ge 1$. Thus, by Theorem A.1 in \cite{AAM:non-autonomous conjugation}, there exists a bounded sequence $\seq{h}$ of holomorphic germs i.e., there exists $\tilde{r}>0$ and holomorphic maps $h_n:B(0;\tilde{r}) \to \C^k$ satisfying $\norm{h_n(z)}<M \norm{z}$ on $B(0;\tilde{r})$ such that 
\[h_{n+1} \circ \w{F}_n= \w{H}_n \circ h_n.\]
Further, we may assume $h_n(0)=0$, $Dh_n(0)=\textsf{Identity}$. By Lemma 5.2 in \cite{AAM:non-autonomous conjugation}, the abstract basins of attraction of $\seq{\w{F}}$ and $\seq{\w{H}}$ at the origin are biholomorphic. Now as in the case $k=2$, it follows that the basins of attraction of $\seq{\w{F}}$ and $\seq{\w{H}}$ at the origin are biholomorphic. Since $\basin{\seq{\w{F}}}=\basin{\seq{S}}$ and $\basin{\seq{\w{H}}}=\basin{\seq{\s}}$, we have $\Omega_{\seq{S}}\cong \Omega_{\seq{\s}}$. Now the proof follows from the next proposition.
\end{proof}

  \begin{prop}\label{p:na p-weak-shift basin}
  Let $\seq{S}$ be a perturbed sequence of weak shift-like maps of sufficiently large degree $d\ge \tilde{d}+2$ associated to the given sequence $\seq{\mathsf{S}}$. Then the basin of attraction of $\seq{{S}}$ at the origin
 \[ \Omega_{\seq{{S}}}=\{z \in \C^k:{S}(n)(z) \to 0 \text{ as }n \to \infty\}\]
 is biholomorphic to $\C^k$.
  \end{prop}

 The proof of this proposition is similar to that of Theorem 1.1 from \cite{B:semigroup}, but the situation here is more involved. Let us first recall a few definitions:
   \begin{align*}
\Omega_{S(m)}=\{z \in \C^k: S(m)^n (z) \to 0 \text{ as } n \to \infty\}, \;
\Omega_{\seq{S}}=\{z \in \C^k: S(n) (z) \to 0 \text{ as } n \to \infty\}.
\end{align*} 
For every $m\ge 1$, note that the origin is an attracting fixed point for $S(m)$. Recall the notation
$s_m=S(m)$ from the previous proof in Section \ref{s:4}. Let $\Omega_{s_m}=\Omega_{S(m)}$ denote the basin of attraction of $s_m$ at the origin, which is a Fatou-Bieberbach domain. Further, it is completely invariant under $S(m)$, i.e.,
\begin{align}\label{e:invariance}
S(m)(\Omega_{s_m})=s_m(\Omega_{s_m})=\Omega_{s_m}=S(m)^{-1}(\Omega_{s_m})=s_m^{-1}(\Omega_{s_m}).	
\end{align}

 \begin{lem}\label{l:compacts in basin}
Let $K$ be a compact set contained in $\Omega_{\seq{S}}$. There exists a positive integer $N_0 = N_0(K) \ge 1$ such that $K \subset \Omega_{s_m}$ for every $m \ge N_0(K).$
\end{lem}
\begin{proof}
  Note that there exists $0<\tilde{r}<\mathsf{r}$ and constants $0<A<\A$ and $0<\B<B<1$, where $\mathsf{r}, \A,\B$ are as in the assumption of Theorem \ref{t:na weak-shift basin}, such that 
  \[A \supnorm{z} \le |\pi_k \circ S_n(z)| \le B\supnorm{z} \]
  for every $z \in B(0;\tilde{r})$. Hence $B(0;\tilde{r})$ is contained in the attracting basin of the origin for every $S_n$, $n \ge 1$. Thus, $\zball{\tilde{r}}\subset \Omega_{s_m} \subset K_m^+\subset\mathcal{K}_m^+$ for every $m \ge 1$. Since $K \subset \Omega_{\seq{S}}$ is compact, there exists $m \ge N_0(K)$ such that $S(m)(w) \in B(0;\tilde{r})$ for every $w \in K$. Hence $S(m)(K) =s_m(K)\subset \Omega_{s_m}$, i.e., by (\ref{e:invariance}) $K \subset \Omega_{s_m}$ for every $m \ge N_0(K).$
\end{proof}

 \begin{proof}[Proof of Proposition \ref{p:na p-weak-shift basin}] 
 Let $\seq{K}$ be an exhaustion by compacts of $\Omega_{\seq{S}}.$ Then from Lemma \ref{l:compacts in basin}, there exists an increasing sequence of positive integers $\seq{m}$ such that $K_n \subset \Omega_{S(m_n)}.$ Note that  every $\Omega_{S(m_n)}$ is a Fatou-Bieberbach domain. Our goal here is to construct a sequence of biholomorphisms $\seq{\phi}$ from $\Omega_{S(m_n)}$ to $\C^k$, appropriately and inductively, such that for a given summable sequence of positive real numbers $\seq{\rho}$ the following holds:
 	\begin{align}\label{e:main}
 		 \norm{\phi_n(z)-\phi_{n+1}(z)}< \rho_n \text{ for } z \in K_n \text{ and }
 		\norm{\phi_n^{-1}(z)-\phi_{n+1}^{-1}(z)}< \rho_n \text{ for } z \in B(0;n).
 	\end{align}
 	
 \noindent	\textit{Basic step: }Let $\phi_1: \Omega_{S(m_1)} \to \C^k$ be a biholomorphism. By results in \cite[Theorem 2.1]{AL:inventiones}, for $\delta< \rho_1/2$ there exists $F_2 \in \Aut$ such that 
 	\begin{align}\label{e:B1}
 	 \norm{\phi_1(z)-F_{2}(z)} &< \delta \text{ for } z \in K_1(r) \text{ and } 
 		\norm{\phi_1^{-1}(z)-F_{2}^{-1}(z)} &< \delta \text{ for } z \in B(0;1+r)
 	\end{align}
 	where $K_1(r)=\cup_{z \in K_1} (B(z;r)) \subset \Omega_{S(m_1)}$ for some $r>0$, i.e., $K_1(r)$ is an $r$-neighbourhood of $K_1$ contained in $\Omega_{S(m_1)}$. Since $F_2$ is uniformly continuous on $K_1(r)$, there exists an $\epsilon_0>0$ such that for $z,w \in K_1(r)$
 	\begin{align}\label{e:B2}
 		\norm{F_2(z)-F_2(w)}< \delta \text{ whenever }\norm{z-w}< \epsilon_0.
 	\end{align}
 	Let $\epsilon< \min\{\epsilon_0,r, \delta\}.$ From \cite[Lemma 4]{W:FB domains}, there exists a biholomorphism $\psi_2: \Omega_{S(m_2)} \to \C^k$ such that
 	\begin{align}\label{e:B3}
  \norm{\psi_2(z)-z}< \epsilon \text{ for every } z \in K_1 \text{ and }
 		\norm{\psi_2^{-1}(z)-z}< \epsilon \text{ for every } z \in F_2^{-1}(B(0;1)).
 	\end{align}
 	Thus, for $z \in K_1$, $\psi_2(z) \in K_1(r)$ and by (\ref{e:B2}), (\ref{e:B3}) it follows that $\norm{F_2 \circ \psi_2(z)-F_2(z)}< \delta.$ Hence from (\ref{e:B1}), 
 	\[\norm{F_2\circ \psi_2(z)-\phi_1(z)}<2\delta< \rho_1 \text{ for }z \in K_1.\]
 	Also by (\ref{e:B3}), for $z \in B(0;1)$, $\norm{\psi_2^{-1} \circ F_2^{-1}-F_2^{-1}(z)}< \epsilon < \delta$. Again by (\ref{e:B1}), 
 	\[\norm{\psi_2^{-1}\circ F_2^{-1}(z)-\phi_1^{-1}(z)}<2\delta< \rho_1 \text{ for }z \in B(0;1).\]
 	Thus $\phi_1$ and $\phi_2:=F_2 \circ \psi_2$ satisfies (\ref{e:main}) for $n=1.$
 	
 	\medskip\noindent
 	\textit{Induction step: }Let $N \ge 2$ and suppose that for all $1 \le n \le N$ there exist biholomorphisms $\phi_n: \Omega_{S(m_n)} \to \C^k$ that satisfy (\ref{e:main}). The goal is to construct $\phi_{N+1}$ such that (\ref{e:main}) holds for $n=N.$ 
 	
 	\medskip
 	
 	As before, for $\delta< \rho_{N}/2$, there exists $F_{N+1} \in \Aut$ such that 
 	\begin{align*}
 		\norm{\phi_{N}(z)-F_{N+1}(z)}< \delta \text{ for } z \in K_{N}(r) \text{ and } 
 		\norm{\phi_{N}^{-1}(z)-F_{N+1}^{-1}(z)}&< \delta \text{ for } z \in B(0;N+r),
 	\end{align*}
 	where $K_{N}(r)$ is an $r$-neighbourhood of $K_{N}$, contained in $\Omega_{S(m_N)}$  for some $r>0$. Since $F_{N+1}$ is uniformly continuous on $K_{N}(r)$, there exists $\epsilon_0>0$ such that for $z,w \in K_N(r)$
 	\begin{align}\label{e:G2}
 		\norm{F_{N+1}(z)-F_{N+1}(w)}< \delta \text{ whenever }\norm{z-w}< \epsilon_0.
 	\end{align}
 	Let $\epsilon< \min\{\epsilon_0,r, \delta\}.$ Then again by \cite[Lemma 4]{W:FB domains}, there exists a biholomorphism $\psi_{N+1}: \Omega_{S(m_{N+1})} \to \C^k$ such that
 	\begin{align}\label{e:G3}
	 \norm{\psi_{N+1}(z)-z}< \epsilon \text{ for } z \in K_{N} \text{ and }
	\norm{\psi_{N+1}^{-1}(z)-z}< \epsilon \text{ for } z \in F_{N+1}^{-1}(B(0;N)). 	
	\end{align}
 	Thus, for $z \in K_N$, $\psi_{N+1}(z) \in K_N(r)$ and by (\ref{e:G2}), (\ref{e:G3}) it follows that $$\norm{F_{N+1} \circ \psi_{N+1}(z)-F_{N+1}(z)}< \delta.$$ Hence 
 	$\displaystyle\norm{F_{N+1}\circ \psi_{N+1}(z)-\phi_N(z)}<2\delta< \rho_N \text{ for }z \in K_N.$
 	Also, similarly as above, by (\ref{e:G3}), for $z \in B(0;N)$, $\displaystyle\norm{\psi_{N+1}^{-1} \circ F_{N+1}^{-1}(z)-F_{N+1}^{-1}(z)}< \epsilon < \delta$ and  by assumption on $F_{N+1}$, 
	\[\norm{\psi_{N+1}^{-1}\circ F_{N+1}^{-1}(z)-\phi_N^{-1}(z)}<2\delta< \rho_N.\]
	Thus $\phi_N$ and $\phi_{N+1}:=F_{N+1} \circ \psi_{N+1}$ satisfies (\ref{e:main}) for $n=N.$

 \medskip As $\seq{\rho}$ is summable, the sequences $\seq{\phi}$ and $\seq{\phi^{-1}}$ constructed above converge on every compact subset of $\Omega_{\seq{S}}$ and $\C^k$. Let the limit maps be $\phi: \Omega_{\seq{S}} \to \C^k$ and $\tilde{\phi}: \C^k \to \C^k$. Being a limit of biholomorphisms, $\phi $ is either injective degenerate in the sense that $\det D\phi \equiv 0$.
 	
 \medskip Choose $M>0$ and choose an $n \ge 1$ sufficiently large, so that $\sum_{i=n}^\infty \rho_i<M/2$. Also, let $K=\phi_n^{-1}(B(0;M))$. Then $\V{K}>0$ and 
\[B(0;M/2) \subset \phi(K) \subset B(0;3M/2)\]
by (\ref{e:main}). Therefore, $\V{\phi(K)}> \V{B(0;M/2)}.$ But if $\det D\phi\equiv 0$, then $\V{\phi(K)}=0$, which is not true. Hence $\phi$ is injective on $\Omega_{\seq{S}}.$ Now by the following Lemma \ref{l:phi onto} it follows that $\Omega_{\seq{S}}$ is biholomorphic to $\C^k$, which completes the proof.
\end{proof}
	
\begin{lem}\label{l:phi onto}
The map $\phi : \Omega_{\seq{S}} \rightarrow \C^k$ is surjective.
\end{lem}

\begin{proof}
 First note that as a consequence of Theorem 5.2 in \cite{DE:Michael Problem}, $\phi_n^{-1}$ converges uniformly to $\tilde{\phi}$ on compact subsets of $\C^k$ and $\tilde{\phi}^{-1}=\phi$ on $\Omega_{\seq{S}}.$  Also $\Omega_{\seq{S}} \subset \tilde{\phi}(\C^k)$.
 	 	 
 \medskip\noindent\textit{Step 1: }For every positive integer $N_0 \ge 1$, $\tilde{\phi}(B(0,N_0)) \subset \mathcal{K}^+$. 
 	 
 \medskip Suppose not. Then there exists $z_0 \in \tilde{\phi}(B(0,N_0))$ such that $z_0 \notin \mathcal{K}^+.$ Let $\tilde{w}_0=\tilde{\phi}^{-1}(z_0) \in B(0;N_0)$ and $z_m=\phi_m^{-1}(\tilde{w}_0).$ Then $z_m \in \Omega_{S(m)}=\Omega_{s_m}$, $\G_m(z_m)=0$, for every $m \ge 1$ and $z_m \to z_0$ by (\ref{e:main}). But by Theorem \ref{t:na weak-shift basin}, $\G_m$ converges uniformly to $G$ on compact subsets of $\C^k$. Hence for a given $\ep>0$ there exists $m_\ep \ge 1$, sufficiently large such that
   \[|G(z_l)-\G_m(z_l)|<\ep \text{ for every } m \ge m_\ep \text{ and }l \ge 1.\] 
  In particular $|G(z_m)|<\ep/2$ for $m \ge m_\ep$. Now as $z_m \to z_0$ and $G$ is continuous on $\C^k$ it follows that $|G(z_0)| \le \ep$. But $\ep>0$ was arbitrary hence $G(z_0)=0$. Thus by (\ref{e:Green}), $z_0 \in \mathcal{K}^+$.
 
 \medskip The above observation also shows that $z_0 \notin \partial \mathcal{K}^+$, for if not, there would exist $\tilde{z}_0 \in \tilde{\phi}(B(0;N_0))$ such that $\tilde{z_0}\in U^+$, and this would contradict the fact that $\tilde{\phi}(B(0;N_0)) \subset \mathcal{K}^+.$ 
  
  \medskip Our next goal is to refine the observation of \textit{Step 1}, i.e.,
  
  \medskip\noindent\textit{Step 2: }For every positive integer $N_0 \ge 1$, $\tilde{\phi}(B(0,N_0)) \subset \I{K^+}.$

  \medskip
  We will prove this by contradiction. So suppose $\tilde{\phi}(B(0;N_0)) \cap \left(\I{\mathcal{K}^+}\setminus {K^+} \right) \neq \emptyset$.  Let 
  \[z_0 \in \tilde{\phi}(B(0;N_0)) \cap \left(\I{\mathcal{K}^+}\setminus {K^+} \right)\] and $\{z_m\}$ be as defined above, a sequence converging to the point $z_0$ with $z_m \in \Omega_{S(m)}$.
 	 	
\medskip\noindent\textit{Claim: }There exists a subsequence $\seq{m}$ such that $\norm{S(m_n)(z_{m_n})} \to \infty$ as $n \to \infty.$

\medskip Choose $\delta_0>0$ such that $B(z_0;\delta_0) \subset \I{\mathcal{K}^+}\setminus {K^+}.$ Since $ z_0\in  \I{\mathcal{K}^+}\setminus {K^+}$, there exists a subsequence $\seq{m}$ such that $S(m_n)(z_0) \to \infty$. As observed in Lemma \ref{l:compacts in basin}, there exists $0<\tilde{r}<\mathsf{r}$ such that $B(0;\tilde{r})\subset \Omega_{\seq{S}}.$ Since $B(z_0;\delta_0) \subset \C^k \setminus K^+$, $S(m)(z) \notin B(0;\tilde{r})$ for every $z \in B(z_0;\delta_0)$ and $m \ge 1$. In particular, $\norm{S(m_n)(z)}>\tilde{r}$ for every $z \in B(z_0;\delta_0)$.

\medskip Thus, $\log \left(\norm{S(m_n)(z)}/{\tilde{r}}\right)$ is a sequence of positive harmonic functions on $B(z_0;\delta_0)$. As $\norm{S(m_n)(z_0)}$ diverges to infinity, Harnack's inequalities show that $\norm{S(m_n)(z)} \to \infty$ uniformly on $B(z_0;\delta_0)$. This proves the claim.

\medskip Choose $\w{R}>\max\{\supnorm{z_n},(k-1)^3,R: n \ge 1\}$ sufficiently large, where $R>1$ is as obtained in Lemma \ref{l:filtration radius} and Remark \ref{r:modified radius}. By the above \textit{claim}, there exists a large $n \ge 1$ such that $\supnorm{S(m_n)(z_{m_n})}>\w{R}$. Suppose $s_{m_n}(z_{m_n})=S(m_n)(z_{m_n}) \in V_{\w{R}}^+ \cup V_{\w{R}}^-.$ Then either 
\[z_{m_n} \in U_{m_n}^+\text{ or }z_{m_n} \in V_{\w{R}^2}^-,
\text{ i.e., } z_{m_n} \notin K_{m_n}^+ \text{ or }\supnorm{z_{m_n}}\ge \w{R}^2,\] respectively. Note that neither of the above conclusions are true by the assumptions on the sequence $\seq{z}$. Hence $S(m_n)(z_{m_n}) \notin V_{\w{R}}^+ \cup V_{\w{R}}^-$ for large enough $n \ge 1$.

\medskip\noindent Fix one such point, say $w_0=S(m_{n_0})(z_{m_{n_0}})$ where $n_0$ is large. 
Further let $m_{n_0}=\en_0$. 
Now $w_0 \in \Omega_{s_{\en_0}}\subset K_{\en_0}^+$ and hence $s_{\en_0}^l(w_0) \to 0$ as $l \to \infty.$ Recall that 
\[
s^{\ku_0+1}_{\en_0} = s_{\en_0} \circ s^{\ku_0}_{\en_0} = S_{\en_0} \circ \cdots \circ S_1 \circ s_{\en_0}^{\ku_0}.
\] 
and the non-autonomous sequence $\seq{\mathbf{S}^{\en_0}}$ defined as 
\[\mathbf{S}^{\en_0}_n= S_{n \text{ mod }\en_0}.\]
Also for positive integers $m,l \ge 1$ and $z \in \C^k$, the notations
\[\mathbf{S}^{\en_0}(m)(z)=\mathbf{S}^{\en_0}_m \circ \cdots \circ \mathbf{S}^{\en_0}_1(z) \text{ and } \mathbf{S}^{\en_0}(m,l)(z)=\mathbf{S}^{\en_0}_{l+m} \circ \cdots \circ \mathbf{S}^{\en_0}_{l+1}(z).\]
Since origin is an uniformly attracting fixed point for the sequence $\seq{\mathbf{S}^{\en_0}}$, there exists a largest positive integer $\ku_0 \ge 0$ such that
\begin{align}\label{e:en_0 bound}
   \supnorm{\mathbf{S}^{\en_0}(\ku_0)(w_0)}\ge \w{R} \text{ and }\supnorm{\mathbf{S}^{\en_0}(m)(w_0)}<\w{R} 
\end{align}
for every $m\ge \ku_0+1.$ 
Let $\we=\mathbf{S}^{\en_0}(\ku_0)(w_0)$ and $\we_i=\pi_i({\we})$ for every $1 \le i \le k$. Note 
\begin{align*}
\mathbf{S}^{\en_0}_n\co=\big(z_2,\hdots,z_{k-1},z_k+Q_{d-1}(z_2),\alpha_nz_1+\tilde{p}_n(z_2,\hdots,z_k)+\ho_{d}(z_2,\hdots, z_k)\big)	
\end{align*}
where 
\[\alpha_n=a_{n \text{ mod } \en_0}\text{ and }\tilde{p}_n(z_2,\hdots,z_k)=p_{n \text{ mod } \en_0}(z_2,\hdots,z_k)\]
with $a_n$ and $p_n$ as defined in (\ref{e:weak shifts}). Hence for every $1 \le i \le k-2$
\[|{\we}_{i+1}|=|\pi_{i}(\mathbf{S}^{\en_0}_{\ku_0+1}({\we}))|=|\pi_i(\mathbf{S}^{\en_0}(\ku_0+1)(w_0))|<\w{R}.\] 
As $\supnorm{\we}\ge \w{R}$ from above either $|\we_1|\ge \w{R}$ or $|\we_k|\ge \w{R}$ or both. If $|{\we}_1|\ge |{\we}_k|$, then $(k-1)|{\we}_1|> |{\we}_k|$. In particular, $\we \in W_{\w{R}}^-$, as defined in Remark \ref{r:modified radius}. Also by the same it follows that $\supnorm{z_{m_{n_0}}}>\w{R}$, which is a contradiction to the choice of $\w{R}$. Hence $|\we_k|=\supnorm{\we}\ge \w{R}$. Also, as $\we \notin V_{\w{R}}^+ \cup  V_{\w{R}}^- $, there exists $1 \le i_0 \le k-1$ such that 
\[|\we_k|<(k-1)|{\we}_{i_0}|.\]

\noindent\textit{Case 1: }Suppose $i_0=1$. Then as discussed above, by Remark \ref{r:modified radius}, $\we \in W_{\w{R}}^-$, we have $\supnorm{z_{m_{n_0}}}>\w{R}$, which is a contradiction.

\medskip\noindent\textit{Case 2: }Suppose $i_0=2.$ Then $\we_{i_0}=\pi_2 (\we )$. Now as $\w{R}>(k-1)^3$ and $d\ge 3$
\begin{align*}
\supnorm{\mathbf{S}^{\en_0}(\ku_0+1)(w_0)}&\ge |\pi_{k-1}(\mathbf{S}^{\en_0}(\ku_0+1)(w_0))|=|\pi_{k-1}(\mathbf{S}^{\en_0}_{\ku_0+1})(\we)|=|\we_k-Q_{d-1}(\we_{i_0})| \\&\ge \bigg|\bigg(\frac{|\we_k|}{k-1}\bigg)^{d-1}-|\we_k| \bigg| \ge |\we_k|\bigg|\bigg(\frac{\w{R}}{k-1}\bigg)^{d-2}-1\bigg)\bigg| 
\\ &>|\we_k|\big((k-1)^{2d-4}-1\big)>\w{R}.
\end{align*}
This contradicts (\ref{e:en_0 bound}).

\medskip\noindent\textit{Case 3: }Suppose $2 < i_0 <k$. Then note that
\[\we_{i_0}=\pi_2 \big(\mathbf{S}^{\en_0}(\ku_0+i_0-2,\ku_0)(\we )\big)=\pi_2\big(\mathbf{S}^{\en_0}(\ku_0+i_0-2)(w_0)\big).\]
Let $\omega_k=\pi_k\big(\mathbf{S}^{\en_0}(\ku_0+i_0-2)(w_0)\big)$. Then by (\ref{e:en_0 bound})
\[
|\omega_k|=|\pi_k\big(\mathbf{S}^{\en_0}(\ku_0+i_0-2)(w_0)\big)| \le \supnorm{\mathbf{S}^{\en_0}(\ku_0+i_0-2)(w_0)}<\w{R} \le |\we_k|.
\]
Now as before, since $\w{R}>(k-1)^3$, $d\ge 3$ and $|\omega_k|<\w{R}\le |\we_k|$,
\begin{align*}
\supnorm{\mathbf{S}^{\en_0}(\ku_0+i_0-1)(w_0)}&\ge \big|\pi_{k-1}\circ \mathbf{S}^{\en_0}_{\ku_0+i_0-1}\circ \mathbf{S}^{\en_0}(\ku_0+i_0-2)(w_0)\big|\\&= \big|\pi_{k-1}\circ \mathbf{S}^{\en_0}_{\ku_0+i_0-1}\circ \mathbf{S}^{\en_0}(\ku_0+i_0-2,\ku_0)(\we)\big|\\&=|\omega_k-Q_{d-1}(\we_{i_0})| \ge  \bigg|\bigg(\frac{|\we_k|}{k-1}\bigg)^{d-1}-|\omega_k|\bigg|\\&\ge \bigg|\bigg(\frac{|\we_k|}{k-1}\bigg)^{d-1}-|\we_k| \bigg| \ge |\we_k|\bigg|\bigg(\frac{\w{R}}{k-1}\bigg)^{d-2}-1\bigg)\bigg| 
\\ &>|\we_k|\big((k-1)^{2d-4}-1\big)>\w{R}.
\end{align*}
 Again, this contradicts (\ref{e:en_0 bound}).

\medskip\noindent
Hence our assumption is not true, i.e., $z_0 \notin \I{\mathcal{K}^+}\setminus {K^+}.$ This proves \textit{Step 2}.

\medskip Finally, if $z_0 \in \partial{K}^+$, then as before there exists $\tilde{z}_0 \in \tilde\phi(B(0;N_0))$ such that $\tilde{z}_0 \in \I{\mathcal{K}^+} \setminus {K}^+.$ Thus, $\tilde{\phi}(B(0;N_0)) \subset \I{K^+}$ for every positive integer $N_0 \ge 1$.

\medskip Now if $z_0 \in \partial \Omega_{\seq{S}} \cap \I{K^+}$, then there exists a neighbourhood of $B(z_0;\delta) \subset \I{K^+}$ where the sequence $\{S(n)\}$ is locally uniformly bounded and hence normal. Then there exists a subsequence $\{l_n\}$ such that $S(l_n)(z_0)$ does not tend to $0$. However, $S(l_n)(z_1)\to 0$  whenever $z_1\in B(z_0;\delta) \cap \Omega_{\seq{S}}$. So, $\{S(n)\}$ is not normal on any neighbourhood of $z_0$. Thus $\partial \Omega_{\seq{S}} \cap \I{K^+}=\emptyset$ and $\tilde{\phi}(\C^k)$ is contained in a single component of $\I{K^+}$. In particular, $\tilde{\phi}(\C^k) \subset \Omega_{\seq{S}}$. But $\phi \circ \tilde{\phi}=\text{Identity}$ on $\C^k$, hence $\tilde{\phi}(\C^k) = \Omega_{\seq{S}}$.
 \end{proof}
\begin{rem}\label{r:Henon}
 Observe that both weak shift-like maps and perturbed weak shift-like maps coincide with the H\'{e}non maps when restricted to two variables. Hence the proof of Theorem 4.1 and Proposition 4.2 also assures that $\Omega_{\h_n} \simeq \C^2$ where $\h_n$'s are uniformly attracting sequence of H\'{e}non maps of the form 
 \[\h_n(x,y)=(y,a_n+p_n(y)),\]
 satisfying the assumptions (i) and (ii), as stated in Section \ref{s:3}. Finally, also note that the assumptions (i) and (ii) implies that the sequence $\{\h_n\}$ satisfies the \emph{uniform filtration and bound condition}, stated in Section \ref{s:2}.
 \end{rem}
 
 \begin{rem*}
Note that in the set up of iterative dynamics of a single H\'{e}non map, the boundary of {\it any} attracting basin is equal to the boundary of the filled Julia set (see \cite[Theorem 2]{BS2}). In the above proof we only show that $\partial \Omega_{\seq{S}}$ is contained in $\partial K^+$. This leads to the question: Is $\partial K^+ = \partial \Omega_{\seq{S}}$ for non-autonomous families of weak shift-like maps or H\'{e}non maps?
 \end{rem*}
 
 \section{Proof of Theorem \ref{t:main result} for \texorpdfstring{$k\ge 3$}{k>2}}\label{s:5}

 Let $\seq{f}$ be a sequence of automorphisms of $\C^k$, $k \ge 3$, admitting a common fixed point at the origin which is uniformly attracting, i.e., there exist constants $A, B$ such that $0<A<B<1$ and $r>0$ such that
 \[ A\norm{z} \le \norm{f_n(z)} \le B\norm{z}\]
for $z \in B(0;r)$. Recall that $k_0 \ge 2$ is the minimum positive integer such that $B^{k_0}<A$. Also, by appealing to the $QR$-factorization of the linear part of $f_n$'s we may assume that $Df_n(0)$ are lower triangular for every $n \ge 1$. In particular, let $Df_n(0)=[u_{ij}^n]$ and $Df_n^{-1}(0)=Df_n(0)^{-1}=[\bar{u}_{ij}^n]$, be the sequences of lower triangular matrices and $\seq{\mathbf{u}}$ and $\seq{\tilde{\mathbf{u}}}$ be the sequence of diagonal elements of $Df_n(0)$ and $Df_n(0)^{-1}$, respectively, i.e.,  
\[
\mathbf{u}_n=({u_{11}^n},\hdots,{u_{kk}^n}), \tilde{\mathbf{u}}_n=(\bar{u}_{11}^n, \hdots,\bar{u}_{kk}^n )=({1}/{u_{11}^n},\hdots, {1}/{u_{kk}^n})
\]
both of which form a sequence of vectors in $\C^k$. Note that $A\le |u_{ii}^n|\le B$ and $0\le |u_{ij}^n|\le B$ for every $1 \le j \le i \le k$.

\medskip We will use the notation $\ze_i=(z:z_i)=(z_1,\hdots,z_{i-1},z_{i+1},\hdots z_k) \in \C^{k-1}$ for every $z \in \C^k$ and $1 \le i \le k$. For a holomorphic polynomial $P:\C^{k-1} \to \C$ and $a=(a_1,\hdots,a_k)\in \C^{k}$, define the following automorphisms of $\C^k$ 
\[T^i_{P,a}(z)=(z_1,\hdots,z_{i-1},a_i z_i+P(\ze_i), z_{i+1},\hdots,z_k).\]
for every $1 \le i \le k$.
Also, we will need to suitably order the set of $k$-tuples of positive integers. For this, recall the lexicographic ordering (from the right side) on $\mathbb{Z}_+^k$. For a pair of $k$-tuples of positive integers $(m_1, m_2, \ldots, m_k)$ and
$(l_1, l_2, \ldots, l_k)$, the notation
\[
(m_1, m_2, \ldots, m_k) < (l_1, l_2, \ldots, l_k)
\]
will mean that that there exists a $j_0$ with $1 \le j_0 < k$ such that
$m_{j_0}<l_{j_0}$ and $m_j=l_j$ for $j_0< j \le k$. In our situation, we will need to order them in the reverse direction, i.e., on a finite set of $k$-tuples the largest element corresponding to the above lexicographic ordering will be the smallest element. This can be defined by saying that 
\begin{equation}\label{e:lexicography}
(m_1,\hdots,m_k) < (l_1,\hdots,l_k)
\end{equation}
if there exists a $j_0$ with $1 \le j_0 \le k$ such that 
$m_{j_0}>l_{j_0}$ and $m_j=l_j$ for $j_0<j \le  k$. Henceforth, this will be called the {\it reverse lexicographic ordering}.

\begin{prop}\label{p:equivalence 1}
 There exist sequences of polynomial maps $\seq{P^i}$, $1 \le i \le k$, where each $P^i_n : \C^{k-1} \ra \C$ has degree at most $k_0$ and such that the sequence $\seq{g}$ defined as
 \begin{align}\label{e:polynomials 1}
 g_n(z) = T^{k}_{P_n^k,\mathbf{u}_n}\circ \cdots \circ T^{1}_{P_n^1,\mathbf{u}_n}(z)	
\end{align}
 satisfies the following:
 \begin{itemize}
 	\item [(i)] $Dg_n(0)=Df_n(0)$.
 	\item [(ii)] The basins of attraction of $\seq{g}$ and $\seq{f}$ at the origin are biholomorphic, i.e., $\Omega_{\seq{g}}\backsimeq \Omega_{\seq{f}}$.
 \end{itemize}
 \end{prop}
 
\begin{proof}
Let $\mathrm{I}^d_l$  denote the set of all positive (non-linear and linear) indices of degree at most $d$, $d \ge 2$ in $\C^l$, $l \ge 2$. In particular,
\begin{align*}
  \mathrm{I}_k^{k_0} &= \left\{(m_1,\hdots, m_k)\in \mathbb{N}^k: 1 \le m_1 + m_2 + \ldots + m_k\le k_0 \right\}, \\
  \mathrm{I}_{k-1}^{k_0} &= \left\{(m_1,\hdots, m_{k-1})\in \mathbb{N}^{k-1}: 1 \le m_1 + m_2 + \ldots + m_{k-1}\le k_0 \right\}, \; \text{and}\\
  \mathrm{I}_k^{k_0-1}&= \left\{(m_1,\hdots, m_k)\in \mathbb{N}^k: 1 \le m_1 + m_2 + \ldots + m_k \le k_0-1 \right\}.
\end{align*}
For every $n \ge 1$ and $1 \le i \le k$, let  
\begin{align}\label{e:Pn}
    P_n^i(\ze_i)=\sum_{j \in \mathrm{I}_{k-1}^{k_0}}\alpha^i_{j,n}\ze_i^j
\end{align}
and $\seq{X} : \C^k \rightarrow \C^k$, a sequence of maps of the form
\begin{align}\label{e:Xn}
\pi_i\circ X_n\co=z_i\left(1+\sum_{\mathsf{j} \in \mathrm{I}_{k}^{k_0-1}} \rho_{\mathsf{j},n}^iz^\mathsf{j}\right).
\end{align}

\noindent \textit{Claim: }For every $j \in \mathrm{I}_{k-1}^{k_0}$, $\mathsf{j} \in \mathrm{I}^{k_0-1}_{k}$ and $1 \le i \le k$, there exist bounded sequences $\{\alpha_{j,n}^i\}$ and $\{\rho_{\mathsf{j},n}^i\}$ such that the maps $\{X_n\}$ and $\{g_n\}$, as in (\ref{e:Xn}) and (\ref{e:polynomials 1}), satisfy 
\begin{align}\label{e:na relation}
	 [X_{n+1}]_{k_0} =[g_n \circ X_n \circ f_n^{-1}]_{k_0}.
\end{align}	
Note that for a fixed $i$, $1 \le i \le k$, by fixing $m_i=0$, we may assume $\mathrm{I}^{k_0}_{k-1} \subset \mathrm{I}^{k_0}_{k}$. We denote this collection of tuples by $\mathrm{I}^{k_0,i}_{k}$ and note that $\mathrm{I}^{k_0,i}_{k} \subset \mathrm{I}^{k_0}_{k}$. Similarly, for a fixed $i$, $1 \le i \le k$, define $\mathrm{I}_{k,i}^{k_0}$, which can be regarded as a subset of $\mathrm{I}^{k_0}_{k}$, as the collection of tuples $m=(m_1, \ldots, m_i, \ldots, m_k) \in \mathrm{I}^{k_0}_k$ with $|m|\ge 2$ such that $(m_1, \ldots, m_i -1, \ldots, m_k) \in \mathrm{I}^{k_0-1}_k$. In particular,
\begin{align*}
 \mathrm{I}^{k_0,i}_{k}&=\{(m_1,\hdots,m_i,\hdots,m_k) : m_i=0\text{ and } 1 \le m_1+m_2+\cdots+m_k\le k_0\}\\   
 \mathrm{I}^{k_0}_{k,i}&=\{(m_1,\hdots,m_i,\hdots,m_k) : m_i\ge 1 \text{ and } 2 \le m_1+m_2+\cdots+m_k\le k_0\}.
\end{align*}  
Let $\mathbf{e}_i=(0,\hdots,1,\hdots,0)$, $1 \le i \le k$ denote the standard basis elements of $\C^k$. Hence 
\begin{align*} 
\mathrm{I}_k^1=\cup_{i=1}^k \mathbf{e}_i \text{ and } \mathrm{I}^{k_0,i}_{k} \cup \mathrm{I}^{k_0}_{k,i}\cup \mathbf{e}_i=\mathrm{I}^{k_0}_{k}.
\end{align*}
Also note the relations (\ref{e:Pn}) and (\ref{e:Xn}) can thus be re-written as follows
\[   P_n^i(z)=\sum_{j \in \mathrm{I}_{k}^{k_0,i}}\alpha^i_{j,n}z^j \text{ and } \pi_i\circ X_n\co=z_i +\sum_{\mathsf{j} \in \mathrm{I}_{k,i}^{k_0}} \rho_{\mathsf{j},n}^iz^\mathsf{j}.\]
Further, the above identities hold for any $k_0 \ge 2$. Since $k_0\ge 2$ is fixed by our set up, we may assume that for any $2 \le l \le k_0$
\begin{align*} 
\left( \mathrm{I}^{l,i}_{k} \setminus \mathrm{I}^{l-1}_k\right) \cup \left(\mathrm{I}^{l}_{k,i} \setminus \mathrm{I}^{l-1}_k\right)=\mathrm{I}^{l}_{k} \setminus \mathrm{I}^{l-1}_{k} \text{ with } \mathrm{I}^{l,i}_{k} \cap \mathrm{I}^{l}_{k,i}=\emptyset
\end{align*}
Let $\sharp \mathrm{I}^{l}_{k}= N_l$ and let the sequences $\seq{\Al^{i,l}}$, $\seq{\Rh^{i,l}}$, for every $1 \le i \le k$ and $1 \le l \le k_0$, in $\C^{N_l}$, be defined as 
\[\rho_{j,n}^i=\begin{cases}
0 &\text{ for } j \in \mathrm{I}^{l,i}_{k}\\
\rho_{j,n}^i &\text{ for } j \in \mathrm{I}_{k,i}^{l}\\
1 &\text{ for }j=\mathbf{e}_i
\end{cases}
\text{ and } \alpha_{j,n}^i=\begin{cases} \alpha_{j,n}^i &\text{ for } j \in \mathrm{I}^{l,i}_{k}\\
 0	&\text{ for } j \in \mathrm{I}_{k,i}^{l}\\
u_{ii}^{n} &\text{ for }j=\mathbf{e}_i
\end{cases}
,\]
where $\Rh_n^{i,l}=(\rho_{j,n}^i)_{j \in \mathrm{I}^{l}_k}$ and $\Al_n^{i,l}=(\alpha_{j,n}^i)_{j \in \mathrm{I}^{l}_k}$. Thus, our goal is to construct bounded sequences $\seq{\Rh^{i,k_0}}$ and $\seq{\Al^{i,k_0}}$, for every $1 \le i \le k$, in $\C^{N_{k_0}}$. Recall that $0<A<B<1$,
\[
\tilde{\mathbf{u}}_n=({1}/{u_{11}^n},\hdots, {1}/{u_{kk}^n}) \in \C^k 
\]
and $A\le |u_{ii}^n|\le B$. For a $k$-tuple of positive integers $m = (m_1, m_2, \ldots, m_k)$, let
\begin{align}\label{e:bar u_n^m}
\tilde{\mathbf{u}}^m_n=\left(\frac{1}{u_{11}^n}\right)^{m_1}\left(\frac{1}{u_{22}^n}\right)^{m_2}\cdots \left(\frac{1}{u_{kk}^n}\right)^{m_k} \in \C. 
\end{align}

Hence, there exists a constant $M>1$ such that 
\[1<B^{-1}\le |\tilde{\mathbf{u}}_n^m|=\left| \left(\frac{1}{u_{11}^n}\right)^{m_1}\left(\frac{1}{u_{22}^n}\right)^{m_2}\cdots \left(\frac{1}{u_{kk}^n}\right)^{m_k}\right |<M \] for every $m=(m_1,\hdots,m_k) \in \mathrm{I}_{k}^{k_0}$ and $n \ge 1$.
Also, $\mathrm{I}_k^l \subset \mathrm{I}_k^{l+1}$ for $1 \le l< k_0$. For a multi-index $m=(m_1,\hdots,m_k)$ with $|m|=1$, i.e., $m=\mathbf{e}_j$, $1 \le j \le k$, we already know $\rho_{\mathbf{e}_j,n}^i=\delta_{ij}.$ Further, $\alpha_{\mathbf{e}_i,n}^i=u_{ii}^n$, by definition, and $\alpha_{\mathbf{e}_j,n}^i=0$ for $i<j\le k$, as $Dg_n(0)$ is a lower triangular matrix. Thus, the sequences $\{\alpha^1_{m,n}\}$ are well-defined for $|m|=1$. Now let
{\small \[\Me_{n}^i=DT_{P_n^i,\mathbf{u}_n}^i(0)=\begin{pmatrix}\mathbf{e}_1\\\vdots \\ \mathbf{e}_{i-1}\\  \Al_n^{i,1}\\ \mathbf{e}_{i+1}\\ \vdots \\ \mathbf{e}_{k}\end{pmatrix} \text{ and } U_n^i=
\begin{pmatrix} u_{11}^n &0 &\cdots &0&0 &0 &\cdots &0\\ \vdots &\vdots & &\vdots &\vdots &\vdots & &\vdots\\ u_{(i-1)1}^n &u_{(i-1)2}^n &\cdots &u_{(i-1)(i-1)}^n &0 &0 &\cdots &0\\ u_{i1}^n &u_{i2}^n &\cdots &u_{i(i-1)}^n&u_{ii}^n &0 &\cdots &0\\ & &&\mathbf{e}_{i+1}\\ &&&\vdots\\ &&&\mathbf{e}_{k}\end{pmatrix} .\]
}
Since $\Me_n^1$ is a well-defined sequence from above, define the sequence of vectors $\Al_n^{i,1}$, inductively such that
\begin{align}\label{e:linear_identity}
\
\Me_n^i \circ \cdots \circ \Me_n^1=U_n^i	
\end{align}
for every $1 \le i \le k$. Then, $\{\Me_n^i\}$ and $\seq{U^i}$ are sequences of invertible matrices, since $\pi_i(\Al_n^{i,1})=u_{ii}^n \neq 0$. Thus, there exist bounded sequences $\{\Al^{i,1}_n\}$ for every $1 \le i \le k$. Further,
\begin{align}\label{e:linear}
Dg_n(0)=\Me_n^k \circ \cdots \circ \Me_n^1=U_n^k=[u_{ij}]=Df_n(0).
\end{align}
 For a map $F$, let $\hat{F}$ denote its non-linear part and $\hat{F}^d$ its $d$-th degree part at the origin. Then, we have
\begin{align*}
X_n(z)=\left(\textsf{Id}+\hat{X}_n\right)(z), \; f_n=\left(Df_n(0)+\hat{f}_n\right)(z), \text{ and }g_n=\Big(Dg_n(0)+\hat{g}_n\Big)(z).
\end{align*}
Thus, from (\ref{e:linear}), the equality of the $k_0$-jets as in (\ref{e:na relation}) reduces to 
\begin{align*}
 \hat{X}_{n+1}\left(Df_n(0)+\hat{f}_n\right)(z)&=Dg_n(0)\circ \hat{X}_n(z)+\hat{g}_n\left(\textsf{Id}+\hat{X}_n\right)(z)-\hat{f}_n(z) \\
 &=Df_n(0)\circ \hat{X}_n(z)+\hat{g}_n\left(\textsf{Id}+\hat{X}_n\right)(z)-\hat{f}_n(z).
\end{align*}
 For $2\le d\le k_0$
\begin{align}\label{e:degree d}
\hat{X}^d_{n+1}(z)=\left(Df_n(0)\circ \hat{X}_n^d\circ Df_n^{-1}(0)\right)(z)+\left(\hat{g}_n^d \circ Df_n^{-1}(0)\right)(z)+\q_n^d(z)
 \end{align}
where $\q_n^d$ is a polynomial map of homogeneous degree $d$ that is determined completely by the sequences $\seq{\Al^{i,d-1}}$ and $\seq{\Rh^{i,d-1}}$, $1 \le i \le k$ and the coefficients of $f_n$, $f_n^{-1}$ (up to degree $d$). Now, as $Df_n(0)$ is a lower triangular matrix, from (\ref{e:degree d}),  $\pi_i \circ \hat{X}^d_{n+1}(z)$ is determined by the following (along with the coefficients of $f_n$ and $f_n^{-1}$, up to degree $d$):
\begin{itemize}
    \item[(i)] $\pi_i \circ \hat{X}^d_{n}(z)$ and $\pi_i \circ \hat{g}^d_{n}(z)$,
    \item[(ii)] $\pi_j \circ \hat{X}^d_{n}(z)$  where $1 \le j \le i-1$,
    \item[(iii)] $\pi_i \circ \q_n^d(z)$.
\end{itemize}
Also, again, due to the fact $Df_n^{-1}(0)$ is a lower triangular matrix, one needs an enumeration of the monomials of $\pi_i \circ \hat{X}_n^d$ and $\pi_i \circ \hat{g}_n^d$ according to which their respective coefficients are evaluated. In particular, let $\phi_d$ be the enumeration induced by the reverse lexicographic ordering of $k-$tuples on the index sets $\mathcal {I}^{d}_k$ for $2 \le d \le k_0$, where 
\[\mathcal{I}_{k}^d=\mathrm{I}_k^d \setminus \mathrm{I}_k^{d-1}= \left\{ (m_1,\hdots,m_k)\in \mathbb{N}^k: m_1 + m_2 + \ldots + m_k =d \right\}.\]
Further, for every $1 \le i \le k$ and $2 \le d \le k_0$
\begin{align}\label{e:disjoint}
\mathcal{I}^d_k\cap \mathrm{I}^d_{k_0,i}=\mathcal{I}^d_k\setminus \mathrm{I}^{d,i}_{k_0}.
\end{align}
Let $g_{n,i}=T^i_{P_n^i,\mathbf{u}_n}\circ \cdots \circ T^1_{P_n^1,\mathbf{u}_n}$, for every $1 \le i\le k$. Then
\begin{align}\label{e:matrix gn}
Dg_{n,i}(0)=D\left (T^i_{P_n^i,\mathbf{u}_n}\circ \cdots \circ T^1_{P_n^1,\mathbf{u}_n}\right)(0)=\begin{pmatrix}
	T_i^n & O_{i\times k-i}\\
	O_{k-i \times i} & \textsf{Id}_{k-i}
\end{pmatrix},	
\end{align}
where $T_i^n$ is a lower triangular matrix for every $n \ge 1$,  and $O$, $\textsf{Id}$ are the zero and identity matrices of appropriate dimensions. Note
\[\pi_{j} \circ g_{n,i}(z)=\pi_{j} \circ g_n (z)\text{ for }1 \le j \le i\]
and 
\[\pi_{j} \circ g_{n,i}(z)=z_j=\pi_j(z) \text{ for }i+1 \le j \le k\]
Hence to determine the sequences $\{\alpha_{m,n}^{i}\}$, $m \in \mathcal{I}_k^d$ it is enough to study the maps $g_{n,i}$, provided we have information about $g_{n,i-1}$, with the assumption that $g_{n,0}=\textsf{Identity}$.

\smallskip
 Note that by (\ref{e:linear}) we know the sequences $\seq{\Al^{i,1}}$ and $\seq{\Rh^{i,1}}$. Let $2 \le d \le k_0$. For $1 \le i \le k $, let $w=g_{n,i-1}(z)$ and $w_j=\pi_j \circ g_{n,i-1}(z)$. Then from above $w_j=z_j$ for every $i \le j \le k$. Now
\[\pi_{i} \circ g_{n,i}(z)=\pi_{i} \circ T^{i}_{P_n^{i},\un} \circ g_{n,i-1}(z).\]
In particular, 
\[[\pi_{i} \circ g_{n,i}(z)]_d=u_{ii}^n z_{i}+ \left[\sum_{m \in \mathrm{I}^{d-1}_k} \alpha ^{i}_{m,n} w^m+\sum_{m \in \mathcal{I}^d_k} \alpha ^{i}_{m,n} w^m\right]_d,\]
or equivalently,
\begin{align}\label{e:g hat}
\pi_{i} \circ \hat{g}^d_{n}(z)= \sum_{m \in \mathcal{I}^d_k} \alpha ^{i}_{m,n} \left(Dg_{n,i-1}(0)z\right)^m+ q_{n,i}^d(z).
\end{align}
where $q_{n,i}^{d}$ is a homogeneous polynomial of degree $d$ depending on the coefficients $f_n$, $f_n^{-1}$ (up to degree d) and the sequences $\seq{\Al^{j,d-1}}, \seq{\Rh^{j,d-1}},\{\alpha^1_{p,n}\},\{\rho^1_{p,n}\},\hdots,\{\alpha^{i-1}_{p,n}\}$ and $\{\rho^{i-1}_{p,n}\}$ for every $1 \le j \le k$ and $p \in \mathcal{I}_k^d$. Thus from (\ref{e:degree d}) and (\ref{e:g hat}) it follows that 
{\small \begin{align}\label{e:degree d_final}
\nonumber \pi_i \circ \hat{X}^d_{n+1}(z)&=\pi_i \left(Df_n(0)\circ \hat{X}_n^d\circ Df_n^{-1}(0)\right)(z)+\left(\pi_i \circ \hat{g}_n^d \circ Df_n^{-1}(0)\right)(z)+\pi_i \circ \q_n^d(z)\\
&=\pi_i \left(Df_n(0)\circ \hat{X}_n^d\circ Df_n^{-1}(0)\right)(z)+\sum_{m \in \mathcal{I}^d_k} \alpha ^{i}_{m,n} \left(Dg_{n,i-1}(0).Df_n^{-1}(0)z\right)^m+\tilde{q}_{n,i}^d(z)
\end{align}
}
where (as before) $\tilde{q}_{n,i}^{d}$ is a homogeneous polynomial of degree $d$ depending on the coefficients $f_n$, $f_n^{-1}$ (up to degree d) and the sequences $\seq{\Al^{j,d-1}}, \seq{\Rh^{j,d-1}},\{\alpha^1_{p,n}\},\{\rho^1_{p,n}\},\hdots,\{\alpha^{i-1}_{p,n}\}$ and $\{\rho^{i-1}_{p,n}\}$ for every $1 \le j \le k$ and $p \in \mathcal{I}_k^d$. 

\smallskip Since $Dg_{n,i-1}(0).Df_n^{-1}(0)$ are lower triangular matrices for every $1 \le i \le k$, we will need an order on the sequence coefficients $\{\alpha_{p,n}^{i}\}$ and $\{\rho_{p,n}^{i}\},p \in \mathcal{I}^d_k$, in which they are evaluated. Again, we will use the enumeration $\phi_d$ induced by the lexicographic ordering (\ref{e:lexicography}) on the set $\mathcal{I}^d_k$. Note that $\sharp \mathcal{I}^d_k=n_d=N_d-N_{d-1}$ and
\begin{align}\label{e:order phi}
\phi_d(1)=(0,\cdots,0,d),\phi_d(2)=(0,\cdots,1,d-1), \hdots, \phi_d(n_d)=(d,0,\cdots,0).
\end{align}
\begin{lem}\label{l:i=i,m=d}
Let $1 \le i \le k$ and $2 \le d \le k_0$. Suppose $m \in \mathcal{I}^d_k$ is such that $\phi_{d}(l_0)=m$, $1 \le l_0\le n_d$. Then
\begin{align}\label{e:main relation}
\rho_{m,n+1}^{i}=u_{ii}^n\tilde{\mathbf{u}}_n^m \rho_{m,n}^{i} +\tilde{\mathbf{u}}_n^m (u_{11}^n)^{m_1} \ldots (u_{(i-1)(i-1)}^n)^{m_i}\alpha_{m,n}^{i}+\li_{n,m,d}^{i},
\end{align}
 where $\tilde{\mathbf{u}}_n^m$ is as defined in (\ref{e:bar u_n^m}) and $\li_{n,m,d}^{i}$ is a bounded term  depending on the following:
 \begin{itemize}
     \item [(i)] the coefficients of $f_n$ and $f_n^{-1}$ up to degree $d$,
     \item[(ii)] the sequences $\seq{\Al^{j,d-1}}$ and $\seq{\Rh^{j,d-1}}$, $1 \le j \le k$, 
     \item[(iii)] the sequences $\{\alpha^\je_{p,n}\},\{\rho^\je_{p,n}\}$ where $1 \le \je \le i-1$ and $p \in \mathcal{I}^d_k$,
     \item[(iv)] the sequences $\{\alpha^{i}_{\phi(l),n}\},\{\rho^{i}_{\phi(l),n}\}$ for every $1 \le l < l_0$.
 \end{itemize} 
 \end{lem} 
 
 \begin{proof}
     Note for $i=1$, $\li_{n,m,d}^i$ is determined by the sequences listed in (i), (ii) and (iv) above. Recall $Df_n(0)=[u_{ij}^n]$, a lower triangular matrix for every $n \ge 1$. Hence
  {\small   \begin{align*}
       \pi_i \left(Df_n(0)\circ \hat{X}_n^d\circ Df_n^{-1}(0)\right)(z)&=\sum_{\je=1}^i u_{i\je}^n \left(\pi_\je\circ \hat{X}_n^d \circ Df_n^{-1}(0)(z)\right)  \\
       &=\sum_{\je=1}^{i-1} u_{i\je}^n \left(\pi_\je\circ \hat{X}_n^d \circ Df_n^{-1}(0)(z)\right)+u_{ii}^n \left(\pi_i\circ \hat{X}_n^d \circ Df_n^{-1}(0)(z)\right)    
     \end{align*}
     }
     Observe that $u_{i\je}^n \left(\pi_\je\circ \hat{X}_n^d \circ Df_n^{-1}(0)(z)\right)$, $1 \le \je \le i-1$ is determined by the coefficients of $f_n$, $f_n^{-1}$ and the sequences $\{\rho^\je_{p,n}\}$ where $1 \le \je \le i-1$ and $p \in \mathcal{I}^d_k$. Thus 
     \begin{align*}
         \pi_i \left(Df_n(0)\circ \hat{X}_n^d\circ Df_n^{-1}(0)\right)(z)=\q_{X,n}^d(z)+\sum_{m \in \mathcal{I}_k^d} u_{ii}^n\rho^i_{m,n} \left(Df_n^{-1}(0)(z)\right)^m
     \end{align*}
     where $\q_{X,n}^d$ is a homogeneous polynomial of degree $d$ that is determined completely by the sequences $\{\rho^\je_{p,n}\}$ where $1 \le \je \le i-1$ and $p \in \mathcal{I}^d_k$. Now from (\ref{e:degree d_final})
     \begin{align*}
     \pi_i \circ \hat{X}^d_{n+1}(z)= \sum_{p \in \mathcal{I}^k_d} \rho^i_{p,n+1} z^{p}&=\sum_{p \in \mathcal{I}^k_d} u_{ii}^n\rho^i_{p,n} \left(Df_n^{-1}(0)(z)\right)^p+\q_{X,n}^d(z) \\ 
      &+\sum_{p \in \mathcal{I}^d_k} \alpha ^{i}_{p,n} \left(Dg_{n,i-1}(0).Df_n^{-1}(0)z\right)^p+\tilde{q}_{n,i}^d(z).
     \end{align*}
     Recall from (\ref{e:linear_identity}) and (\ref{e:linear}) that 
     \[Df_n(0)=U_n^k\text{ and }Dg_{n,i-1}(0)=U_n^{i-1}.\] Let $A^{i-1}_n=(U^{i-1}_n)^{-1}U^k_n$. Then the above reduces to
     \begin{align}\label{e:phi_d sum}
   \sum_{l=1}^{n_d}\rho^i_{\phi_d(l),n+1} (z)^{\phi(l)}&=\sum_{l=1}^{n_d}\left\{u_{ii}^n\rho^i_{\phi_d(l),n}  \big(U_n^k(z)\big)^{\phi_d(l)} +\alpha ^{i}_{\phi_d(l),n} \left(A^{i-1}_n (z)\right)^{\phi_d(l)}\right\}+\tilde{\q}_{n}^d(z)  
   \end{align}
where $\tilde{\q}_{n}^d(z)$ is a homogeneous polynomial of degree $d$ depending on the sequences (i)--(iii) in the statement of the lemma.

\medskip Now let $1 \le l_1 < l_2 \le n_d$ and assume that $(r_1,\hdots,r_k)=\phi_d(l_1)<\phi_d(l_2)=(s_1,\hdots,s_k)$, as per the ordering defined in (\ref{e:lexicography}). In particular, there exists a $j_0$ with $1 \le j_0 \le k $ such that 
\[r_j=s_j \text{ for }j_0<j\le k, \text{ and } r_{j_0}>s_{j_0}.\]
For a lower triangular matrix $T=[t_{ij}]$, define
\begin{align}\label{e:triangular z}
  (T(z))^{\phi_d(l_2)}&=(t_{11}z_1)^{s_1}(t_{21}z_1+t_{22}z_2)^{s_2}\hdots(t_{k1}z_1+\cdots +t_{kk}z_k)^{s_k}=\sum_{l \ge l_2} \tau_{\phi_d(l)} z^{\phi_d(l)}.
\end{align}

\noindent {\it Claim}: The coefficient of the monomial $z^{\phi_d(l_1)} = z_1^{r_1}z_2^{r_2} \ldots z_k^{r_k}$ in $(T(z))^{\phi_d(l_2)}$ is zero. 

\medskip

 Suppose $p=(p_1,\hdots,p_k) \in \mathcal{I}^d_k$ is such that the expression for the coefficient of the monomial $z^p = z_1^{p_1}z_2^{p_2} \ldots z_k^{p_k}$ in (\ref{e:triangular z}) is non-zero. Then we have the following straightforward

\medskip\noindent\textit{Induction statement:} Suppose that $p_{k-j}=s_{k-j}$ for every $0 \le j \le  \je$ and every $\je \ge 0$.  Then 
\[p_{k-\je-1} \le s_{k-\je-1}.\]
Thus, if the coefficient of $z^{\phi_d(l_1)} = z_1^{r_1}z_2^{r_2} \ldots z_k^{r_k}$ is non-zero in (\ref{e:triangular z}) then $r_{j_0} \le s_{j_0}$, which is not true as we have assumed $l_1<l_2$. This proves the claim.

\medskip Also note that by (\ref{e:matrix gn}), the diagonal elements of $U^k_n$ and $A^{i-1}_n$ are
\[\un=(u_{11}^n,\hdots,u_{kk}^n) \;\text{and } \mathbf{a}^{i-1}_n=\left(1,\hdots,1, \frac{1}{u_{ii}^n},\hdots, \frac{1}{u_{kk}^n}\right).\] 
Recall that $\phi_d(l_0)=m=(m_1,\hdots, m_k)$. Thus, the coefficient of $z^{\phi_d(l_0)}=z^m$ in $\left(U_n^k(z)\right)^{m}$ and $\left(A_n^{i-1}(z)\right)^{m}$
are $\un^{m}$ and 
\[(\mathbf{a}_n^{i-1})^m=\tilde{\mathbf{u}}_n^{m}(u_{11}^{n})^{m_1}\hdots (u_{(i-1)(i-1)}^n)^{m_{i-1}},\]
respectively. Now, by the above claim
\begin{align*}
    \left(U_n^k(z)\right)^{m}&= \un^{m}z^m+\sum_{l >l_0}b_n^l z^{\phi_d(l)} \text{ and }\\
    \left(A_n^{i-1}(z)\right)^{m}&= \tilde{\mathbf{u}}_n^{m}(u_{11}^{n})^{m_1}\hdots (u_{(i-1)(i-1)}^n)^{m_{i-1}}z^m+\sum_{l >l_0}c_n^l z^{\phi_d(l)} 
\end{align*}
where $b_n^l$ and $c_n^l$ depend on the linear coefficients of $f_n$ and $f_n^{-1}$ for every $l_0<l\le n_d$. Thus, from (\ref{e:phi_d sum}) we have 
\[\rho_{m,n+1}^{i}=u_{ii}^n\tilde{\mathbf{u}}_n^m \rho_{m,n}^{i} +\tilde{\mathbf{u}}_n^m (u_{11}^n)^{m_1}...(u_{(i-1)(i-1)}^n)^{m_{i-1}}\alpha_{m,n}^{i}+\li_{n,m,d}^{i},\]
where $\li_{n,m,d}^{i}$ is a bounded term  depending on the sequences (i)--(iv).
 \end{proof}
Finally, we proceed to complete the proof of the proposition.

\medskip\noindent\textbf{Initial Case:} Suppose $d=|m|=2$, i.e., $m\in\mathcal{I}_k^2= \mathrm{I}_k^2\setminus \mathrm{I}_k^1$.

\smallskip Our goal is to construct the sequences $\{\Al^{i,2}_n\}$ and $\{\Rh^{i,2}_n\}$ inductively for every $1 \le i \le k$ in the ascending order of $i$.

\medskip\textit{Basic Step:} Let $i=1$ and fix an $m=(m_1,\hdots,m_k) \in \mathcal{I}_k^2$. From Lemma \ref{l:i=i,m=d} we have the following
\begin{align}\label{e:i=1,d=2}
\rho_{m,n+1}^1=u_{11}^n\tilde{\mathbf{u}}_n^m \rho_{m,n}^1+\tilde{\mathbf{u}}_n^m\alpha_{m,n}^1+\li_{n,m,2}^1,
\end{align}
 where $\phi_2(l_0)=m$, $1 \le l_0\le n_2$ and $\li_{n,m,2}^1$ is a bounded term depending on the following
 \begin{itemize}
     \item [(i)] the coefficients of $f_n$ and $f_n^{-1}$ up to degree $2$,
     \item[(ii)] the sequences $\seq{\Al^{j,1}}$ and $\seq{\Rh^{j,1}}$, $1 \le j \le k$ constructed above, and
     \item[(iii)] the sequences $\{\alpha^1_{\phi(l),n}\},\{\rho^1_{\phi(l),n}\}$ for every $1 \le l < l_0$.
 \end{itemize} 
 
By (\ref{e:disjoint}), for every $m \in \mathcal{I}_k^2$ either $m \in \mathrm{I}^{2,2}_k$ or $m \in \mathrm{I}^2_{k,2}$. As $\rho_{m,n}^1=0$, there exists a bounded sequence $\{\alpha_{m,n}^1\}$ for every $m\in  \mathrm{I}_k^{2,1} \cap \mathcal{I}_k^2$ such that $\rho_{m,n+1}^1=0$. Further 
\[1<|u_{11}^n\tilde{\mathbf{u}}_n^m|<M\] 
for every $m\in  \mathrm{I}_{k,1}^2 \cap \mathcal{I}_k^2$. Hence, by Lemma \ref{l:linear bound}, there exist appropriate values of $\rho_{m,1}^1 \in \C$, for every $m\in  \mathrm{I}_{k,1}^2 \cap \mathcal{I}_k^2$, such that the sequence $\{\rho_{m,n}^1\}$ as defined inductively, is bounded with $\alpha_{m,n}^1=0$ for every $n \ge 1$. Thus the sequences $\{\alpha^1_{\phi_2(l),n}\}$ and $\{\rho^1_{\phi_2(l),n}\}$ can be uniquely determined from (\ref{e:i=1,d=2}) when evaluated in the ascending order $\phi_2(l)$, $1 \le l \le n_2$.

 \medskip Now suppose for a given $1\le i \le k$, there exist bounded sequences $\{\rho_{m,n}^{\mathsf{j}}\}$ and $\{\alpha_{m,n}^{\mathsf{j}}\}$ for every $1 \le \mathsf{j} \le i$ and $m \in  \mathcal{I}_k^2$. Then we have the

\medskip\textit{Inductive Step for $i+1$:} Fix an $m=(m_1,\hdots,m_k)=\phi_2(l_0)\in \mathcal{I}_k^2$, $1 \le l_0\le n_2$. Then as in the previous cases, by Lemma \ref{l:i=i,m=d} we have
\begin{align}\label{e:i=i+1,m=2}
\rho_{m,n+1}^{i+1}=u_{(i+1)(i+1)}^n\tilde{\mathbf{u}}_n^m \rho_{m,n}^{i+1} +\tilde{\mathbf{u}}_n^m (u_{11}^n)^{m_1}...(u_{ii}^n)^{m_i}\alpha_{m,n}^{i+1}+\li_{n,m,2}^{i+1},
\end{align}
 where $\li_{n,m,2}^{i+1}$ is a bounded term depending on the following
 \begin{itemize}
     \item [(i)] the coefficients of $f_n$ and $f_n^{-1}$ up to degree $2$,
     \item[(ii)] the sequences $\seq{\Al^{j,1}}$ and $\seq{\Rh^{j,1}}$, $1 \le j \le k$, 
     \item[(iii)] the sequences $\{\alpha^\je_{p,n}\},\{\rho^\je_{p,n}\}$ where $1 \le \je \le i$ and $p \in \mathcal{I}^2_k$,
     \item[(iv)] the sequences $\{\alpha^{i+1}_{\phi(l),n}\},\{\rho^{i+1}_{\phi(l),n}\}$ for every $1 \le l < l_0$.
\end{itemize}

Thus, arguing similarly as before, that is, evaluating (\ref{e:i=i+1,m=2}) in the ascending order of $m=\phi_2(l) \in \mathcal{I}^2_k$, $1 \le l \le n_2$, there exist sequences $\{\rho_{m,n}^{i+1}\}$ and $\{\alpha_{m,n}^{i+1}\}$ in $\C$ for every $m \in \mathcal{I}^2_k$, defined inductively by the above relations, that are either bounded or identically zero sequences.

\medskip\noindent  Hence by induction,  there exist bounded sequences ${\Al_{n}^{j,2}}$ and ${\Rh_{n}^{j,2}} \in \C^{N_2}$, for every $1 \le j \le k$. 

\medskip Next, we assume there exist bounded sequences ${\Al_{n}^{j,d}}$ and ${\Rh_{n}^{j,d}} \in \C^{N_d}$, for some $2\le d < k_0$ and for every $1 \le j \le k$. We proceed towards the 

\medskip\noindent \textbf{General Case:} Suppose $|m|=d+1$, i.e., $m\in \mathcal{I}^{d+1}_k=\mathrm{I}_k^{d+1}\setminus \mathrm{I}_k^{d}$.
 
 \medskip Note, the \textit{steps} here are exactly similar to those in $|m|=2$. By Lemma \ref{l:i=i,m=d} we work with the identity (\ref{e:main relation}) first in the ascending order of $i$, $1 \le i \le k$. Further for every fixed $i$, we work with $\{\alpha^i_{m,n}\}$ and $\{\rho^i_{m,n}\}$ in the ascending order of $m=\phi_{d+1}(l)$, $1 \le l \le n_{d+1}$. Hence there exist bounded sequences $\{\Al_{n}^{i,d+1}\}$ and $\{\Rh_{n}^{i,d+1}\} \in \C^{N_{d+1}}$ for every $1 \le i \le k$. 
 
 \smallskip Now, by induction on $d$, there exist bounded sequences $\{\Al_{n}^{i,k_0}\}$ and $\{\Rh_{n}^{i,k_0}\} \in \C^{N_{k_0}}$ for every $1 \le i \le k$, which proves the claim.

\medskip Thus, the sequences $\seq{f}$ and $\seq{g}$ are non-autonomously conjugated up to order $k_0$ at the origin. Hence, the proof of the proposition follows as indicated in Section \ref{s:2}.
\end{proof}

\noindent Now we complete the
 \begin{proof}[Proof of Theorem \ref{t:main result} in $\C^k$, $k \ge 3$] Recall the maps
 \[T^i_{P_n^i,\mathbf{u}^n}(z)=(z_1,\hdots,z_{i-1},u_{ii}^n z_i+P^i_n(\ze_i), z_{i+1},\hdots,z_k)\]
 for every $1 \le i \le k$, where $\ze_i=(z_1,\hdots,z_{i-1},z_{i+1},\hdots z_k) \in \C^{k-1}$ and $P_n^i: \C^{k-1} \ra \C$ are polynomial maps of degree at most $k_0$. By Proposition \ref{p:equivalence 1}, given the sequence of automorphisms $\seq{f}$ of $\C^k$, satisfying (\ref{e:ub}), there exists a sequence of polynomial maps $\seq{g}$ from $\C^k$ to $\C^k$ of the form
 \[g_n(z)=T^k_{P_n^k,\mathbf{u}^n}\circ \cdots \circ T^1_{P_n^1,\mathbf{u}^n},\] 
 such that $Dg_n(0)(z)=Df_n(0)(z)$ and $\Omega_{\seq{g}} \backsimeq \Omega_{\seq{f}}$. Also, $k_0 \ge 2$ is chosen large enough so that $B^{k_0}<A$ and $Dg_n(0)$ is a lower triangular matrix for every $n \ge 1$. 
 
 \medskip Now, there exist sequences $\seq{m}$ and $\seq{r}$ such that $n=m_nk-r_n$ where $0\le r_n \le k-1$. For $0 \le i \le k-1$, let
\[P_n^{k-i}(\ze_{k-i})=P_n^{k-i}(z_1,\hdots,z_{k-i-1},z_{k-i+1},\hdots z_k)=\p_{nk-i}(z_{(k-i)+1},\hdots,z_{k},z_1\hdots,z_{k-i-1}).\]
Now let $\seq{\s}$ be a sequence of weak shift-like maps  of the form
\begin{align*}
\s_n\co&=\left(z_2,\hdots,z_{k}, u_{(k-r_n)(k-r_n)}^{m_n}z_1+\p_{m_nk-r_n}(\ze_1)\right)\\ &=\left(z_2,\hdots,z_{k}, u_{(k-r_n)(k-r_n)}^{m_n}z_1+P_{m_n}^{k-r_n}(z_{r_n+2},\hdots,z_{k},z_2,\hdots,z_{r_n+1})\right).
\end{align*}

\noindent {\it Claim}: For every $0\le i \le k-1$,
\[\s_{nk-i}\circ \cdots \circ \s_{(n-1)k+1}(z)=\inv_{k-i} \circ g_{n,(k-i)},\] 
where (recall from the proof of Proposition \ref{p:equivalence 1})
\[ g_{n,(k-i)}=T^{k-i}_{P_n^{k-i},\mathbf{u}^n}\circ \cdots \circ T^1_{P_n^1,\mathbf{u}^n},\] 
and $\inv_{k-i}$ is the inversion of the first $(k-i)$-coordinates and the last $i$-coordinates in $\C^k$, i.e.,
\[\inv_{k-i}(z_1,\hdots,z_k)=(z_{(k-i)+1},\hdots,z_k,z_1,\hdots,z_{k-i}).\] 
We will prove this by induction on $i$, in the descending order, starting from $i=k-1$ to $i=0$. For $i=k-1$ and $n \ge 1$, 
\begin{align*}
\s_{(n-1)k+1}\co&=\big(z_2,\hdots,z_k,u_{11}^n z_1 +P_n^1(z_2,\hdots,z_k)\big) \\
&=\inv_1\big(u_{11}^n z_1 +P_n^1(z_2,\hdots,z_k),z_2,\hdots,z_k\big)=\inv_1 \circ g_{n,1}.	
\end{align*}
Assume that the statement is true for some
$i=k-l$, i.e.,
\begin{align*}
\s_{(n-1)k+l}\circ \cdots \circ\s_{(n-1)k+1}\co=\big(z_{l+1},\hdots, z_k, \pi_1 \circ g_{n,l}(z),\cdots, \pi_{l}\circ g_{n,l}(z)\big).	
\end{align*}
Now  for $i=k-(l+1)$,
\begin{align*}
&\pi_k\circ \s_{(n-1)k+(l+1)} \circ \s_{(n-1)k+l}\circ \cdots \circ\s_{(n-1)k+1}\co \\
&= u_{(l+1)(l+1)}^n z_{l+1}+\p_{(n-1)k+l+1}(z_{l+2},\hdots, z_k, \pi_1 \circ g_{n,l}(z),\cdots, \pi_{l}\circ g_{n,l}(z))\\
&= u_{(l+1)(l+1)}^n z_l+P_n^{l+1}( \pi_1 \circ g_{n,l}(z),\cdots, \pi_{l}\circ g_{n,l}(z),z_{l+2},\hdots, z_k)\\
&= \pi_{l+1}\circ g_{n,(l+1)}\co.
\end{align*}
Note that 
\[\pi_j\circ g_{n,(l+1)}\co=\pi_j \circ g_{n,{l}}\co \] for every $1 \le j \le l$ and 
\[\pi_j\circ g_{n,(l+1)}\co=\pi_j \circ g_{n,{l}}\co =z_j\] for every $l+2 \le j \le k$. 
Hence the above observation gives
\begin{align*}
 &\s_{(n-1)k+(l+1)} \circ \s_{(n-1)k+l}\circ \cdots \circ\s_{(n-1)k+1}(z)\\&=\big(z_{l+2},\hdots, z_k, \pi_1 \circ g_{n,(l+1)}(z),\cdots, \pi_{l+1}\circ g_{n,(l+1)}(z)\big)=\inv_{l+1}\circ g_{n,(l+1)}(z),	
\end{align*}
which proves the induction statement and the claim. Hence, for $n \ge 1$
\[\s_{nk} \circ \s_{(n-1)k+l}\circ \cdots \circ\s_{(n-1)k+1}(z)=g_{n, k}(z)=g_n(z).\]
Now $\seq{\s}$ is a sequence of weak shift-like maps satisfying the criteria of Theorem \ref{t:na weak-shift basin} and hence $\Omega_{\seq{\s}}=\Omega_{\seq{g}}$ is biholomorphic to $\C^k$. Thus, $\Omega_{\seq{f}}$ is biholomorphic to $\C^k$.
 \end{proof}
 \begin{rem*}
Note that the proof of Theorem \ref{t:main result} in Section \ref{s:2}, i.e., in $\C^2$ does not emphasize on the maps $T_{P_n^i,\mathbf{u}_n}(z)$, as we directly compute the sequence $\{g_n\}$ of H\'{e}non maps in $\C^2$ (see Proposition \ref{p:step 1}). However observe that from the identity (\ref{e:rseq})
\begin{align*}  
g_n(x,y)=T_{P_n^2,\mathbf{u}_n} \circ T_{P_n^1,\mathbf{u}_n}(x,y)&=\left(u_{11}^n x+P_n^1(y),u_{22}^n y+P_n^2\left(u_{11}^n x+P_n^1(y)\right) \right)
\\&=\left(y, u_{22}^n+P_n^2(y)\right) \circ \left(y, u_{11}^n+P_n^1(y)\right).
\end{align*}
 \end{rem*}
 
 \section{Explanation of the method for \texorpdfstring{$k=3$}{k=3} and \texorpdfstring{$k_0=5$}{k0=5}}\label{s:6}
 
The purpose of this section is to illustrate the various steps that go into the proof of Proposition $5.1$ by working through a specific example.

\medskip

Let $\seq{f}$ be a sequence of automorphisms of $\C^3$ admitting a common fixed point at the origin that is uniformly attracting, i.e., there exist $0<A<B<1$ and $r>0$ such that
 \[ A\norm{z} \le \norm{f_n(z)} \le B\norm{z}\]
for all $\vert z \vert < r$ and $B^{5}<A\le B^4$. Also, $U_n^3=Df_n(0)$ is a lower triangular $3 \times 3$ matrix, i.e.,
\[U_n^3=\begin{pmatrix}
	u_{11}^n &0 &0 \\
	u_{21}^n &u_{22}^n &0\\
	u_{31}^n &u_{32}^n &u_{33}^n 
\end{pmatrix},\]
where $A\le |u_{ii}^n|\le B$ for $1 \le i \le 3$ and $0 \le |u_{ij}^n|\le B$, for $1 \le i \neq j \le 3$. Let $\mathbf{u}_n=(u_{11}^n,u_{22}^n,u_{33}^n)$ and $\tilde{\mathbf{u}}_n=({1}/{u_{11}^n}, {1}/{u_{22}^n}, {1}/{u_{33}^n})$.  Our goal is to construct
\begin{align}
g_n(z):=T^{3}_{P_n^3,\mathbf{u}_n}\circ T^{2}_{P_n^2,\mathbf{u}_n} \circ T^{1}_{P_n^1,\mathbf{u}_n}(z)	
\end{align}	
where 
\begin{align*}
T^1_{P_n^1}(z)&=(u_{11}^n z_1+P_n^1(z_2,z_3), z_2, z_3)\\ T^2_{P_n^2}(z)&=(z_1, u_{22}^n z_2+P_n^2(z_1,z_3), z_3)\\T^3_{P_n^3}(z)&=(z_1, z_2,u_{33}^n z_3+P_n^3(z_1,z_2)),
\end{align*}
and
\begin{align*}
	P_n^1(z_1,z_2,z_3)&=\sum_{(0,m_2,m_3) \in \mathrm{I}^{5,1}_{2}} \alpha^1_{(0,m_2,m_3),n}z_2^{m_2}z_3^{m_3},\\ P_n^2(z_1,z_2,z_3)&=\sum_{(m_1,0,m_3) \in \mathrm{I}^{5,2}_{2}} \alpha^2_{(m_1,0,m_3),n}z_1^{m_1}z_3^{m_3}, \\ P_n^3(z_1,z_2,z_3)&=\sum_{(m_1,m_2,0) \in \mathrm{I}^{5,3}_{3}} \alpha^3_{(m_1,m_2,0),n}z_1^{m_1}z_2^{m_2}.
\end{align*}
Further, we assume the various components of the $X_n$'s to be given by 
\begin{align*}
\pi_1 \circ X_n(z_1,z_2,z_3)&= z_1+ \sum_{(m_1,m_2,m_3) \in \mathrm{I}^{4}_{3,1}} \rho^1_{(m_1,m_2,m_3),n}z_1^{m_1}z_2^{m_2}z_3^{m_3}\\
\pi_2 \circ X_n(z_1,z_2,z_3)&= z_2+ \sum_{(m_1,m_2,m_3) \in \mathrm{I}^{4}_{3,2}} \rho^2_{(m_1,m_2,m_3),n}z_1^{m_1}z_2^{m_2}z_3^{m_3}\\
\pi_3 \circ X_n(z_1,z_2,z_3)&= z_3+ \sum_{(m_1,m_2,m_3) \in \mathrm{I}^{4}_{3,3}} \rho^3_{(m_1,m_2,m_3),n}z_1^{m_1}z_2^{m_2}z_3^{m_3}.
\end{align*}
Let $m=(m_1,m_2,m_3) \in \mathrm{I}_3^5$. For $|m|=m_1+m_2+m_3=1$, it straightaway follows that
 \[\alpha^1_{(1,0,0),n}=u^n_{11},\;\alpha^2_{(0,1,0),n}=u^n_{22},\; \alpha^3_{(0,0,1),n}=u^n_{33}\]
and 
\[\alpha^1_{(0,1,0),n}=0,\; \alpha^1_{(0,0,1),n}=0, \; \alpha^2_{(0,0,1),n}=0.\]
Note that
 \[DT_{P_n^1,\mathbf{u}_n}(0)=\begin{pmatrix}
	\alpha^1_{(1,0,0),n} &0 &0 \\
	0 &1 &0\\
	0 &0 &1
\end{pmatrix} \text{ and }
D(T_{P_n^1,\mathbf{u}_n} \circ T_{P_n^2,\mathbf{u}_n})(0)=\begin{pmatrix}
	u_{11}^n &0 &0 \\
	u_{21}^n &u_{22}^n &0\\
	0 &0 &1
	\end{pmatrix}.
\]
Hence 
\[D T_{P_n^2,\mathbf{u}_n}(0)=\begin{pmatrix}
	1 &0 &0 \\
	\alpha^2_{(1,0,0),n} &\alpha^2_{(0,1,0),n} &0\\
	0 &0 &1
	\end{pmatrix}=\begin{pmatrix}
	u_{11}^n &0 &0 \\
	u_{21}^n &u_{22}^n &0\\
	0 &0 &1
	\end{pmatrix}
\begin{pmatrix}
	u_{11}^n &0 &0 \\
	0 &1 &0\\
	0 &0 &1
	\end{pmatrix}^{-1}.\]
Similarly
{\small \[D T_{P_n^3,\mathbf{u}_n}(0)=\begin{pmatrix}
	1 &0 &0 \\ 
    0 &1 &0\\
	\alpha^3_{(1,0,0),n} &\alpha^3_{(0,1,0),n} &\alpha^3_{(0,0,1),n}
	\end{pmatrix}=\begin{pmatrix}
	u_{11}^n &0 &0 \\
	u_{21}^n &u_{22}^n &0\\
	u_{31}^n &u_{32}^n &u_{33}^n
	\end{pmatrix}
\begin{pmatrix}
	u_{11}^n &0 &0 \\
	u_{21}^n &u_{22}^n &0\\
	0 &0 &1
	\end{pmatrix}^{-1}.\]
	}
Thus, the linear coefficients of $g_n$'s and $X_n$'s which are given by the sequences $\seq{\Al^{i,1}}, \seq{\Rh^{i,1}}\in \C^3$ are determined uniquely for every $1 \le i \le 3$.

\begin{itemize}
\item Let $|m|=2$. The lexicographic ordering $\phi_2$ on $\mathcal{I}^2_3$, as referred in (\ref{e:lexicography}) is 
{\small \begin{align*}
&{\phi}_2(1)=(0,0,2),\; {\phi}_2(2)=(0,1,1), \;{\phi}_2(3)=(1,0,1), \\ &{\phi}_2(4)=(0,2,0),\; {\phi}_2(5)=(1,1,0),\; {\phi}_2(6)=(2,0,0).
\end{align*}}
\item Let $|m|=3$. The lexicographic ordering $\phi_3$ on $\mathcal{I}^3_3$, as referred in (\ref{e:lexicography}) is as follows
{\small \begin{align*}
{\phi}_3(1)&=(0,0,3),\; {\phi}_3(2)=(0,1,2),\; {\phi}_3(3)=(1,0,2),  {\phi}_3(4)=(0,2,1),
\\{\phi}_3(5)&=(1,1,1),\;{\phi}_3(6)=(2,0,1),\; {\phi}_3(7)=(0,3,0),\; {\phi}_3(8)=(1,2,0),\\
{\phi}_3(9)&=(2,1,0),\;{\phi}_3(10)=(3,0,0).
\end{align*}}
\item Let $|m|=4$. The lexicographic ordering $\phi_4$ on $\mathcal{I}^4_3$, as referred in (\ref{e:lexicography}) is as follows
{\small \begin{align*}
{\phi}_4(1)&=(0,0,4),\; {\phi}_4(2)=(0,1,3),\; {\phi}_4(3)=(1,0,3),\; {\phi}_4(4)=(0,2,2),\\ 
{\phi}_4(5)&=(1,1,2),\;{\phi}_4(6)=(2,0,2),\; {\phi}_4(7)=(0,3,1),\; {\phi}_4(8)=(1,2,1),\\
{\phi}_4(9)&=(2,1,1), \;{\phi}_4(10)=(3,0,1), \; {\phi}_4(11)=(0,4,0),\; {\phi}_4(12)=(1,3,0),\\ 
{\phi}_4(13)&=(2,2,0),\; {\phi}_4(14)=(3,1,0),\; {\phi}_4(15)=(4,0,0).
\end{align*}}
\item Let $|m|=5$. The lexicographic ordering $\phi_5$ on $\mathcal{I}^5_3$, as referred in (\ref{e:lexicography}) is as follows
{\small \begin{align*}
{\phi}_5(1)&=(0,0,5), \;{\phi}_5(2)=(0,1,4),\; {\phi}_5(3)=(1,0,4),  \;{\phi}_5(4)=(0,2,3),\\ 
{\phi}_5(5)&=(1,1,3),\;{\phi}_5(6)=(2,0,3),\; {\phi}_5(7)=(0,3,2),\; {\phi}_5(8)=(1,2,2),\\
{\phi}_5(9)&=(2,1,2),\;{\phi}_5(10)=(3,0,2),\;{\phi}_5(11)=(0,4,1),\;{\phi}_5(12)=(1,3,1),\\ 
{\phi}_5(13)&=(2,2,1),\; {\phi}_5(14)=(3,1,1),\; {\phi}_5(15)=(4,0,1), \;{\phi}_5(16)=(0,5,0),\\
{\phi}_5(17)&=(1,4,0),\; {\phi}_5(18)=(2,3,0),\; {\phi}_5(19)=(3,2,0), {\phi}_5(20)=(4,1,0), \\ 
{\phi}_5(21)&=(5,0,0).
\end{align*}}
\end{itemize}
By Lemma \ref{l:i=i,m=d}, for $1 \le i \le 3$ and $m \in \mathcal{I}^d_3$, $2 \le d \le 5$ such that $\phi_{d}(l_0)=m$, $1 \le l_0\le n_d$. Then
\begin{align}\label{e:example}
\rho_{m,n+1}^{i}=u_{ii}^n\tilde{\mathbf{u}}_n^m \rho_{m,n}^{i} +\tilde{\mathbf{u}}_n^m (u_{11}^n)^{m_1}...(u_{(i-1)(i-1)}^n)^{m_i}\alpha_{m,n}^{i}+\li_{n,m,d}^{i},
\end{align}
where $\li_{n,m,d}^{i}$ is a bounded term  depending on the following
 \begin{itemize}
     \item [(i)] coefficients of $f_n$ and $f_n^{-1}$ up to degree $d$,
     \item[(ii)] sequences $\seq{\Al^{j,d-1}}$ and $\seq{\Rh^{j,d-1}}$, $1 \le j \le 3$, 
     \item[(iii)]sequences $\{\alpha^\je_{p,n}\},\{\rho^\je_{p,n}\}$ where $1 \le \je \le i-1$ and $p \in \mathcal{I}^d_3$,
     \item[(iv)] sequences $\{\alpha^{d}_{\phi(l),n}\},\{\rho^{d}_{\phi(l),n}\}$ for every $1 \le l \le l_0$.
 \end{itemize} 

By (\ref{e:disjoint}), note that for any particular $1 \le i \le 3$ and $m \in \mathrm{I}^5_{3}$ with $|m|=d\ge 2$ if one knows the sequence $\{\li_{n,m,d}^i\}$, it is possible to determine the sequence $\{\alpha_{m,n}^i\}$ and $\{\rho_{m,n}^i\}$ uniquely. Thus we solve (\ref{e:example}) in the following order 
\begin{itemize}
    \item[(i)] $i=1$, evaluate  $\{\alpha_{\phi_2(l),n}^1\}$ and $\{\rho_{\phi_2(l),n}^1\}$ in the ascending order $1 \le l \le 6$.
    \item[(ii)]$i=2$, evaluate  $\{\alpha_{\phi_2(l),n}^2\}$ and $\{\rho_{\phi_2(l),n}^2\}$ in the ascending order $1 \le l \le 6$.
     \item[(iii)]$i=3$, evaluate  $\{\alpha_{\phi_2(l),n}^3\}$ and $\{\rho_{\phi_2(l),n}^3\}$ in the ascending order $1 \le l \le 6$.
    \item[(iv)] $i=1$, evaluate  $\{\alpha_{\phi_3(l),n}^1\}$ and $\{\rho_{\phi_3(l),n}^1\}$ in the ascending order $1 \le l \le 10$.
    \item[(v)]$i=2$, evaluate  $\{\alpha_{\phi_3(l),n}^2\}$ and $\{\rho_{\phi_3(l),n}^2\}$ in the ascending order $1 \le l \le 10$.
     \item[(vi)]$i=3$, evaluate  $\{\alpha_{\phi_3(l),n}^3\}$ and $\{\rho_{\phi_3(l),n}^3\}$ in the ascending order $1 \le l \le 10$.
     \item[(vii)] $i=1$, evaluate  $\{\alpha_{\phi_4(l),n}^1\}$ and $\{\rho_{\phi_4(l),n}^1\}$ in the ascending order $1 \le l \le 15$.
    \item[(viii)]$i=2$, evaluate  $\{\alpha_{\phi_4(l),n}^2\}$ and $\{\rho_{\phi_4(l),n}^2\}$ in the ascending order $1 \le l \le 15$.
     \item[(ix)]$i=3$, evaluate  $\{\alpha_{\phi_4(l),n}^3\}$ and $\{\rho_{\phi_4(l),n}^3\}$ in the ascending order $1 \le l \le 15$.
     \item[(x)] $i=1$, evaluate  $\{\alpha_{\phi_5(l),n}^1\}$ and $\{\rho_{\phi_3(l),n}^1\}$ in the ascending order $1 \le l \le 21$.
    \item[(xi)]$i=2$, evaluate  $\{\alpha_{\phi_5(l),n}^2\}$ and $\{\rho_{\phi_3(l),n}^2\}$ in the ascending order $1 \le l \le 21$.
     \item[(xii)]$i=3$, evaluate  $\{\alpha_{\phi_5(l),n}^3\}$ and $\{\rho_{\phi_5(l),n}^3\}$ in the ascending order $1 \le l \le 21$.
     \end{itemize}
Hence all together, there exists a sequence of polynomials $\{P_n^i\}$, for every $1 \le i \le 5$, of degree at most $5$ such that
\[g_n(z)=(z_1, z_2,u_{33}^n z_3+P_n^3(z_1,z_2))\circ (z_1, u_{22}^n z_2+P_n^2(z_1,z_3), z_3)\circ  (u_{11}^n z_1+P_n^1(z_2,z_3), z_2, z_3)\]
for $z=(z_1,z_2,z_3) \in \C^3$. Further, $\seq{g}$ is non-autonomously conjugated to $\seq{f}$, hence $\Omega_{\seq{g}}$ is biholomorphic to $\seq{f}$. Now, let 
\begin{align*}
 	S_{3n}(z_1,z_2,z_3)&=(z_2,z_3,u_{33}^n z_1+\p_{3n}(z_2,z_3))=(z_2,z_3,u_{33}^n z_1+P_{n}^3(z_2,z_3))\\
 	S_{3n-1}(z_1,z_2,z_3)&=(z_2,z_3,u_{22}^n z_1+\p_{3n-1}(z_2,z_3))=(z_2,z_3,u_{22}^n z_1+P_{n}^2(z_3,z_2))\\
 	S_{3n-2}(z_1,z_2,z_3)&=(z_2,z_3,u_{11}^n z_1+\p_{3n-2}(z_2,z_3))=(z_2,z_3,u_{11}^n z_1+P_{n}^1(z_2,z_3))
 \end{align*}
and note that $g_n(z)=S_{3n} \circ S_{3n-1} \circ S_{3n-2}(z)$. Hence $\seq{g}$ is a sequence of weak shift-like maps and by Theorem \ref{t:na weak-shift basin}, $\Omega_{\seq{g}}$ is biholomorphic to $\C^3$.


\bibliographystyle{amsplain}

\end{document}